\documentclass[final,leqno]{article}
\usepackage{epsfig,psfrag,amsmath,amssymb}
\newtheorem{mytheorem}{Theorem}
\newtheorem{mylemma}{Lemma}

\newtheorem{myremark}{Remark}

\newcommand{\N}{{\sf N\hspace*{-1.0ex}\rule{0.15ex}{1.3ex}\hspace*{1.0ex}}}

\newcommand{\BA}{\begin{eqnarray}}
\newcommand{\EA}{\end{eqnarray}}
\newcommand{\BE}{\begin{equation}}
\newcommand{\EE}{\end{equation}}
\newcommand{\ba}{\begin{array}}
\newcommand{\ea}{\end{array}}
\newcommand{\baa}{\begin{eqnarray*}}
\newcommand{\eaa}{\end{eqnarray*}}
\newcommand{\be}{\begin{equation}}
\newcommand{\ee}{\end{equation}}

\def\N{{\cal N}}

\def\u{{\bf u}}

\def\Rtilde{{\widetilde{R}}}
\def\Btilde{{\widetilde{B}}}

\def\Atilde{{\widetilde{A}}}
\def\Wtilde{{\widetilde{W}}}

\def\>{\raisebox{-1ex}{$\; \stackrel{\textstyle >}{\sim } \; $}}
\def\<{\raisebox{-1ex}{$ \; \stackrel{\textstyle <} {\sim } \; $}}
\newcommand{\vvec}[1]{{\mathbf{#1}}}

\newcommand{\EQ}[1]{(\ref{equation:#1})}
\newcommand{\LA}[1]{\ref{lemma:#1}}

\def\endproof{\qquad $\Box$}
\def\beginproof{\indent {\it Proof:~}}

\begin{document}



\title{A unified FETI-DP approach for incompressible Stokes equations}

\author{Xuemin Tu\thanks{Department of Mathematics, University of Kansas, 1460 Jayhawk Blvd, Lawrence, KS 66045-7594,  {\tt xtu@math.ku.edu}, {\tt http://www.math.ku.edu/$\sim$xtu/}.} \and Jing Li\thanks{Department of Mathematical Sciences, Kent State University, Kent, OH 44242, {\tt li@math.kent.edu}, {\tt http://www.math.kent.edu/$\sim$li/}.} }
\date{}
\maketitle

\begin{abstract}
A unified framework of FETI-DP algorithms is proposed for solving the system of linear equations arising from the mixed finite element approximation of incompressible Stokes equations. A distinctive feature of this framework is that it allows using both continuous and discontinuous pressures in the algorithm, while previous FETI-DP methods only apply to discontinuous pressures. A preconditioned conjugate gradient method is used in the algorithm with either a lumped or a Dirichlet preconditioner, and scalable convergence rates are proved. This framework is also used to describe several previously developed FETI-DP algorithms and greatly simplifies their analysis.
Numerical experiments of solving a two-dimensional incompressible Stokes problem demonstrate the performances of the discussed FETI-DP algorithms represented under the same framework.
\end{abstract}



\section{Introduction}
The finite element tearing and interconnecting (FETI) methods were
introduced by Farhat and Roux \cite{FETI1,FETI2,FETI3,FETI4} for
second order elliptic problems.  In these algorithms, Lagrange
multipliers are introduced to enforce the continuity of the solution
across the subdomain interface. The original system of linear
equations is reduced to a symmetric positive semi-definite system for
the Lagrange multipliers, which can be solved by a preconditioned
conjugate gradient method. Both a lumped preconditioner \cite{FETI2}
and a Dirichlet preconditioner \cite{FETI4} have been used in the FETI
methods. Compared with the lumped preconditioner, the Dirichlet
preconditioner is more effective in the reduction of iteration count, but it is also more expensive. Numerical experiments in \cite{FETI4} show that the lumped preconditioner is more efficient for second-order problems while the Dirichlet preconditioner is better for plate and shell problems.

The dual-primal FETI (FETI-DP) method, introduced by Farhat {\em et. al.} \cite{Farhat:2001:FDP, Farhat:2000:SDP}, represents a further development of the FETI methods. In a FETI-DP  algorithm, a few variables from each subdomain are selected as the coarse level primal variables which are shared by neighboring subdomains, while the continuity of the other subdomain interface variables is enforced by using Lagrange multipliers. The reduced system for the Lagrange multipliers becomes symmetric positive definite and is solved by a preconditioned conjugate gradient method using either the lumped or the Dirichlet preconditioner. With an appropriate choice of coarse level variables, the condition number bound of the FETI-DP algorithm has been proved independent of the number of subdomains for both second-order and fourth-order elliptic systems, and in both two and three dimensions, cf.~\cite{Mandel:2000:CDP, kla02}.

The FETI-DP algorithm was first extended to solving incompressible
Stokes equations by Li~\cite{li05}. In addition to the coarse level
primal velocity variables, the subdomain average pressure degrees of
freedom are also selected as the coarse level variables and the
resulting coarse level problem is  symmetric indefinite. The reduced system for the Lagrange multipliers is still symmetric positive definite and a preconditioned conjugate gradient method can be used. Only the Dirichlet preconditioner was studied in \cite{li05}, and it was proved for two-dimensional problems that, under the condition that both the subdomain corner and certain edge-average velocity degrees of freedom are selected as coarse level primal variables,  the condition number bound is independent of the number of subdomains and grows only polylogarithmically with the size of individual subdomain problems.

Recently, Kim, Lee, and Park \cite{kim10} introduced a different FETI-DP formulation for solving the incompressible Stokes problems, where no pressure variables are selected as coarse level primal variables and the resulting coarse level problem is symmetric positive definite. Only the lumped preconditioner was studied in \cite{kim10}, for which the edge-average velocity degrees of freedom are no longer needed in the coarse level problem; in fact, as few coarse level primal variables as for solving positive definite elliptic problems were used and as strong condition number bound was established. Reduction of the coarse level problem size has also been achieved by Dohrmann and Widlund in an overlapping Schwarz type algorithm for solving almost incompressible elasticity, cf. \cite{Doh09, Doh10}. Keeping the size of the coarse problem small is important in large scale computations; a large coarse problem can be a bottleneck and additional efforts in the algorithm are needed to reduce its impact, cf.~\cite{Tu:2004:TLB,Tu:2005:TLB,Tudd16,Klawonn:2005:IFM,Dohrmann:2005:ABP,KimTu,Tu:2011:TLBS}.

All  above mentioned algorithms and their analysis for solving incompressible Stokes problems are only valid for finite element discretization using discontinuous pressures. Discontinuous pressures have also been used in domain decomposition algorithms for other similar type saddle-point problems; see \cite{Kla98, 8pavwid, li05, paulo2003, Doh04, li06, kim06, Tu:2005:BPP, Tu:2005:BPD, LucaOlofStef}.

In \cite{li12}, the authors proposed a new FETI-DP algorithm for solving incompressible Stokes equations, which allows using both discontinuous and continuous pressures in the finite element discretization.
The lumped preconditioner was studied in \cite{li12} and, similar to \cite{kim10}, as few coarse level primal variables as for solving positive definite elliptic problems were used and as strong condition number bound was established.

The purpose of this paper is two-fold. First, the FETI-DP formulation
proposed in \cite{li12} is used as a unified framework to describe
the two previous FETI-DP algorithms studied in \cite{li05} and
\cite{kim10}. It is observed that these two FETI-DP
algorithms can be represented as special cases of using
discontinuous pressures in the new formulation. The condition number bound estimate based on the new formulation also greatly simplifies the analysis in \cite{li05} and \cite{kim10}. Second, a new Dirichlet preconditioner is studied for the FETI-DP algorithm presented under the unified framework using either continuous or discontinuous pressures. The same condition number bound as in \cite{li05} is obtained. Moreover, this new Dirichlet preconditioner
involves solving symmetric positive definite subdomain problems and is
less expensive compared with
the Dirichlet preconditioner used in \cite{li05} where subdomain saddle-point problems need be solved.
To stay focused on the purpose of this paper, the presentation of the
algorithms and their analysis is restricted to the case of solving
two-dimensional problems.

The rest of this paper is organized as follows. The finite element
discretization of the incompressible Stokes equation is introduced in
Section \ref{section:FEM}. A domain decomposition approach is
described in Section~\ref{section:DDM}. A reduced system of
equations is derived in Section~\ref{section:Gmatrix}. Section
\ref{section:techniques} provides some techniques used in the
condition number bound estimate. The lumped and the Dirichlet
preconditioners are studied in Sections~\ref{section:lumped}
and~\ref{section:Dirichlet}, respectively.  At the end, in Section~\ref{section:numerics}, numerical results for solving a two-dimensional incompressible Stokes problem demonstrate the performances of the discussed algorithms and their connections.

\section{Finite element discretization}
\label{section:FEM}

We consider the following incompressible Stokes problem on a
bounded, two-dimensional polygonal domain $\Omega$  with a
Dirichlet boundary condition,
\begin{equation}
\label{equation:Stokes}
\left\{
\begin{array}{rcll}
-\Delta {\bf u} + \nabla p & = & {\bf f}, & \mbox{ in } \Omega \mbox{ , } \\
-\nabla \cdot {\bf u}       & = & 0, & \mbox{ in } \Omega \mbox{ , } \\
{\bf u}                    & = & {\bf u}_{\partial \Omega}, & \mbox{ on } \partial \Omega \mbox{ . }\\
\end{array}\right.
\end{equation}
The boundary velocity ${\bf u}_{\partial \Omega}$ satisfies the compatibility
condition $\int_{\partial \Omega} {\bf u}_{\partial \Omega} \cdot {\bf
  n} = 0$. Without loss of generality, we assume that ${\bf u}_{\partial \Omega} = {\bf 0}$.

The weak solution of \EQ{Stokes} is given by: find $\vvec{u} \in
\left(H^1_0(\Omega)\right)^2 = \{ \vvec{v} \in (H^1(\Omega))^2
\mid \vvec{v} = \vvec{0} \mbox{ on }
\partial \Omega \}$, and $p \in
L^2(\Omega)$, such that,
\begin{equation}
\label{equation:bilinear} \left\{
\begin{array}{lcll}
a(\vvec{u}, \vvec{v}) + b(\vvec{v}, p) & = & (\vvec{f}, \vvec{v}),
& \forall \vvec{v}\in \left(H^1_0(\Omega)\right)^2 , \\ [0.5ex]
b(\vvec{u}, q) & = & 0, & \forall q \in L^2(\Omega) \mbox{ , }
\\
\end{array} \right.
\end{equation}
where
\[
a(\vvec{u}, \vvec{v})= \int_{\Omega} \nabla{\bf u} \cdot \nabla{\bf v}, \quad
b(\vvec{u},q) = -\int_{\Omega} (\nabla \cdot \vvec{u}) q, \quad
(\vvec{f}, \vvec{v}) = \int_{\Omega} \vvec{f} \cdot \vvec{v}.
\]
The solution of \EQ{bilinear} is not unique, with $p$ different up to an additive constant.

\begin{figure}[hbt]
\begin{center}
\includegraphics[scale=0.55]{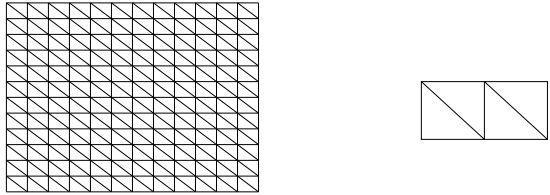}
\end{center}
\caption{A mixed finite element with discontinuous pressures.}
\label{fig:2dfem}
\end{figure}

A mixed finite element method is used to solve \EQ{bilinear}. $\Omega$ is triangulated into shape-regular elements of characteristic size $h$. $\vvec{W} \in \left(H^1_0(\Omega)\right)^2$ represents the velocity finite element space and it contains continuous functions. The pressure finite element space is represented by $Q \subset L^2(\Omega)$.  Both continuous and discontinuous pressures can be used in our algorithm. A mixed finite element space with discontinuous pressure is shown on the left in Figure~\ref{fig:2dfem} on a uniform triangular mesh of a rectangular domain, where the velocity is piecewise linear on the mesh and the pressure is a constant on each union of four triangles as shown on the right in the figure.
This mixed finite element was used in \cite{li05}.

\begin{figure}[h]
\begin{center}
\includegraphics[scale=.65]{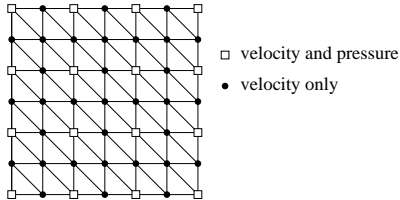}
\end{center}
\caption{\label{figure:TaylorHood} The modified Taylor-Hood mixed finite element.}
\end{figure}

A mixed finite element space with continuous pressure is the modified Taylor-Hood mixed finite element, as shown in Figure~\ref{figure:TaylorHood}, cf.~\cite[Chapter III, \S 7]{braess}. The velocity finite element space contains the piecewise linear functions on the finest triangular mesh and the pressure finite element space contains the piecewise linear functions on the coarser triangular mesh with the doubled mesh size.

The finite element solution $(\vvec{u}, p) \in \vvec{W} \bigoplus Q$ of \EQ{bilinear} satisfies
\begin{equation}
\label{equation:matrix} \left[
\begin{array}{cccc}
A     &  B^T \\
B     &  0   \\
\end{array}
\right] \left[
\begin{array}{c}
{\bf u}     \\
p   \\
\end{array}
\right] = \left[
\begin{array}{l}
{\bf f}        \\
0      \\
\end{array}
\right] .
\end{equation}
Here $A$, $B$ and $\vvec{f}$  represent respectively the restrictions of $a( \cdot , \cdot )$, $b(\cdot, \cdot )$ and $\vvec{f}( \cdot)$ to the finite-dimensional spaces $\vvec{W}$ and $Q$. In this paper the same notation is used to represent both a finite element function and the vector of its nodal values.

The coefficient matrix in \EQ{matrix} is rank deficient. $A$ is symmetric positive definite. The kernel of $B^T$, denoted by $Ker(B^T)$, is the space of all constant pressures in $Q$. The range of $B$, denoted by $Im(B)$,  is orthogonal to $Ker(B^T)$ and is the subspace of $Q$ consisting of all vectors with zero average. The solution of \EQ{matrix} always exists and is uniquely determined when the pressure is considered in the quotient space $Q/Ker(B^T)$. In this paper, when $q \in Q/Ker(B^T)$, $q$ always has zero average. For a more general right-hand side vector $({\bf f}, ~ g)$ given in \EQ{matrix}, the existence of its solution requires that $g \in Im(B)$, i.e., $g$ has zero average.

Both mixed finite elements shown on Figures \ref{fig:2dfem} and
\ref{figure:TaylorHood} are inf-sup stable in the sense that there
exists a positive constant $\beta$, independent of $h$, such that
\begin{equation}
\label{equation:infsup}
\sup_{{\bf w} \in {\bf W}}
\frac{\left< q, B \vvec{w} \right>^2}{\left< \vvec{w},  A \vvec{w} \right>} \geq \beta^2 \left< q, Z q \right>,
\hspace{0.5cm} \forall q \in Q/Ker(B^T).
\end{equation}
Here, as always in this paper, $\left< \cdot, \cdot \right>$ represents the inner product of two vectors. The matrix $Z$ represents the mass matrix defined on the pressure finite element space, i.e., for any $q \in Q$, $\|q\|_{L^2}^2 =
\left< q, Z q \right>$. $Z$ is spectrally equivalent to $h^2 I$, where $I$ represents the identity matrix of the same dimension, i.e., there exist positive constants $c$ and $C$, such that
\be
\label{equation:massmatrix}
c h^2 I \leq Z \leq C h^2 I,
\ee
cf. \cite[Lemma B.31]{Toselli:2004:DDM}. Here, and in other places of this paper, $c$ and $C$ represent generic positive constants which are independent of the mesh size $h$ and $H$ (discussed in the following section).

\section{A non-overlapping domain decomposition approach}
\label{section:DDM}

The domain $\Omega$ is decomposed into $N$ nonoverlapping polygonal subdomains $\Omega_i$, $i = 1, 2, ..., N$. Each subdomain is a union of a bounded number of elements, with the diameter of each subdomain in the order of $H$. The nodes on the boundaries of
neighboring subdomains match across the subdomain interface $\Gamma =
{(\cup\partial\Omega_i)} \backslash
\partial\Omega$. $\Gamma$ is composed of
subdomain edges, which are regarded as open subsets of
$\Gamma$, and of the subdomain vertices, which are end points of
edges.

The velocity and pressure finite element spaces ${\bf W}$ and $Q$ are decomposed into
\[
{\bf W} = {\bf W}_I \bigoplus {\bf W}_{\Gamma},\quad
Q = Q_I \bigoplus Q_\Gamma,
\]
 respectively. Here ${\bf W}_I$ and $Q_I$ are the direct sums of subdomain
interior velocity spaces ${\bf W}^{(i)}_I$, and subdomain interior
pressure spaces $Q^{(i)}_I$, respectively, i.e.,
$$
{\bf W}_I = \bigoplus_{i=1}^{N}{\bf W}^{(i)}_I, \quad
Q_I =
\bigoplus_{i=1}^{N}Q^{(i)}_I.
$$
${\bf W}_{\Gamma}$ is the subdomain boundary velocity space and it contains the subdomain boundary velocity degrees of freedom shared by neighboring subdomains. $Q_\Gamma$ contains the subdomain boundary pressure degrees of freedom shared by neighboring subdomains. For the case of using discontinuous pressures, no pressure degrees of freedom are shared by neighboring subdomains and $Q_\Gamma$ is empty. In fact, for the discontinuous pressure case, the algorithm presented in this paper allows that $Q_\Gamma$ either be empty or contain any number of intrinsically subdomain interior pressure degrees of freedom; more details on this will be discussed in Sections \ref{subsection:cp} and  \ref{subsection:dp}.

To formulate our domain decomposition algorithm, we introduce a partially sub-assembled interface
velocity space 
\[
\vvec{\Wtilde}_{\Gamma} = \vvec{W}_{\Pi} \bigoplus
\vvec{W}_{\Delta} = \vvec{W}_{\Pi} \bigoplus \left(
\bigoplus_{i=1}^N \vvec{W}^{(i)}_\Delta \right).
\]
Here, $\vvec{W}_{\Pi}$ is the continuous coarse level velocity space and the coarse level primal velocity degrees of freedom are shared by neighboring subdomains. The complimentary space $\vvec{W}_{\Delta}$ is the direct sum  of subdomain dual interface velocity spaces $\vvec{W}_{\Delta}^{(i)}$, which correspond to the remaining interface velocity degrees of freedom and are spanned by basis functions which vanish at the primal degrees of freedom. Thus, an element in the space $\vvec{\Wtilde}_{\Gamma}$ has a continuous primal velocity component and typically a discontinuous dual velocity component.

In this paper, two choices of $\vvec{W}_{\Pi}$ are used. In the first,
$\vvec{W}_{\Pi}$ is spanned by all the subdomain corner velocity nodal
basis functions and the coarse level primal variables are only the
subdomain corner velocity variables. In the second, besides all the
subdomain corner velocity nodal basis functions, on each edge
$\Gamma^{ij}$ shared by neighboring subdomains $\Omega_i$ and
$\Omega_j$, $\vvec{W}_{\Pi}$ is also spanned by an edge-average finite
element basis function such that  $\int_{\Gamma^{ij}} \vvec{w}^{(i)}_\Delta \cdot \vvec{n}_{ij} = 0$,
for any $\vvec{w}^{(i)}_\Delta \in \vvec{W}^{(i)}_\Delta$. Here $\vvec{n}_{ij}$ denotes a fixed selection of the normal to $\Gamma^{ij}$. Therefore for the second choice of $\vvec{W}_{\Pi}$, the following divergence free boundary condition
\be
\label{equation:divergencefree}
\int_{\partial \Omega_i} {\bf w}_{\Delta}^{(i)} \cdot {\bf n} = 0
\ee
is satisfied for all  $\vvec{w}^{(i)}_\Delta \in
\vvec{W}^{(i)}_\Delta$. For more details
on choosing coarse level primal edge-average velocity variables to
satisfy the divergence free condition for incompressible Stokes
problems, including for the three-dimensional case, see \cite[Section
7]{li06}. We note that the choice of $\vvec{W}_{\Pi}$ depends on the preconditioner used in the algorithm. The first choice is sufficient for using the lumped preconditioner, but for the Dirichlet preconditioner the second one has to be used; for more detailed discussions, see Sections \ref{section:lumped} and \ref{section:Dirichlet}.

The functions ${\bf w}_{\Delta}$ in ${\bf W}_{\Delta}$ are in general not continuous across $\Gamma$. To enforce their continuity, a boolean matrix $B_\Delta$ is constructed from $\{0,1,-1\}$. Each row of $B_\Delta$ only contains two non-zero entries, $1$ and $-1$, corresponding to the same velocity degree of freedom on each subdomain boundary node, but attributed to two neighboring subdomains. For any ${\bf w}_{\Delta}$ in ${\bf W}_{\Delta}$, each row of $B_\Delta {\bf w}_{\Delta} = 0$ implies that the two degrees of freedom from the two neighboring subdomains be the same. When non-redundant continuity constraints are enforced, $B_\Delta$ has full row rank. We denote the range of $B_\Delta$ applied on ${\bf W}_{\Delta}$ by $\Lambda$, the vector space of the Lagrange multipliers.

In order to define a certain subdomain boundary scaling operator, we introduce a positive scaling factor $\delta^{\dagger}(x)$ for each node $x$ on the subdomain boundary $\Gamma$. Let $\N_x$ be the number of subdomains sharing $x$, and we simply take $\delta^{\dagger}(x) = 1/\N_x$. In applications, these scaling factors will depend on the heat conduction coefficient and the first of the Lam\'{e} parameters for scalar elliptic problems and the equations of linear
elasticity, respectively; see \cite{kla02,kla06}. Given such scaling factors on the subdomain boundary nodes, we can define a scaled operator $B_{\Delta, D}$. We recall that each row of $B_\Delta$ has only two nonzero entries, $1$ and $-1$, corresponding to the same subdomain boundary node $x$. Multiplying each entry by the scaling factor $\delta^{\dagger}(x)$ gives us $B_{\Delta, D}$.

Then solving the original fully assembled linear system~\EQ{matrix} is equivalent to: find
$({\bf u}_I, ~p_I, ~{\bf u}_{\Delta}, ~{\bf u}_{\Pi}, ~p_{\Gamma}, ~\lambda) \in
{\bf W}_I \bigoplus Q_I \bigoplus {\bf W}_{\Delta} \bigoplus {\bf W}_\Pi \bigoplus Q_\Gamma \bigoplus \Lambda$, such that
\be
\label{equation:bigeq}
\left[
\begin{array}{cccccc}
A_{II}      & B_{II}^T       & A_{I \Delta}       & A_{I \Pi}       & B_{\Gamma I}^T     &  0           \\[0.8ex]
B_{II}      & 0              & B_{I \Delta}       & B_{I \Pi}       & 0                  &  0           \\[0.8ex]
A_{\Delta I}& B_{I \Delta} ^T& A_{\Delta\Delta}   & A_{\Delta \Pi}  & B_{\Gamma \Delta}^T&  B_{\Delta}^T\\[0.8ex]
A_{\Pi I}   & B_{I \Pi}^T    & A_{\Pi \Delta}     & A_{\Pi \Pi}     & B_{\Gamma \Pi}^T   &  0           \\[0.8ex]
B_{\Gamma I}& 0              & B_{\Gamma \Delta}  & B_{\Gamma \Pi}  & 0                  &  0           \\[0.8ex]
0           & 0              & B_{\Delta}         & 0               & 0                  &  0
\end{array}
\right]
\left[ \begin{array}{c}
{\bf u}_I        \\[0.8ex]
p_I              \\[0.8ex]
{\bf u}_{\Delta} \\[0.8ex]
{\bf u}_{\Pi}    \\[0.8ex]
p_{\Gamma}     \\[0.8ex]
\lambda
\end{array} \right] =
\left[ \begin{array}{l}
{\bf f}_I        \\[0.8ex]
0                \\[0.8ex]
{\bf f}_{\Delta} \\[0.8ex]
{\bf f}_\Pi      \\[0.8ex]
0                \\[0.8ex]
0
\end{array} \right] \mbox{ .  }
\ee

Corresponding to the one-dimensional null space of~\EQ{matrix}, we consider a vector of the form
$\left( \u_I,~p_I, ~\u_\Delta, ~ \u_\Pi, ~p_\Gamma, ~\lambda \right)  = \left( {\bf 0},~1_{p_I},~{\bf 0},
~{\bf 0}, ~1_{p_\Gamma}, \lambda \right)$, where $1_{p_I}$ and $1_{p_\Gamma}$ represent vectors with value $1$ on each entry of the vector. Substituting it into \EQ{bigeq} gives zero blocks on the  right-hand side, except at the third block
\be
\label{equation:fdelta}
{\bf f}_{\Delta} = [B_{I\Delta}^T ~~ B_{\Gamma\Delta}^T]\left[\begin{array}{c}1_{p_I}\\
1_{p_\Gamma}\end{array}\right]+B^T_\Delta \lambda.
\ee
The first term in \EQ{fdelta} represents the line integral of the normal components of the velocity finite element basis functions across the subdomain boundary, and corresponding to the same velocity degree of freedom on the subdomain boundary, their values on the two neighboring subdomains are negative of each other. Therefore
\[
[B_{I\Delta}^T ~~ B_{\Gamma\Delta}^T]\left[\begin{array}{c}1_{p_I}\\
1_{p_\Gamma}\end{array}\right] = B^T_\Delta B_{\Delta,D}[B_{I\Delta}^T ~~ B_{\Gamma\Delta}^T]\left[\begin{array}{c}1_{p_I}\\
1_{p_\Gamma}\end{array}\right],
\]
from which we know that $\vvec{f}_\Delta = {\bf 0}$ in \EQ{fdelta}, for
\[
\lambda =-B_{\Delta,D}[B_{I\Delta}^T ~~ B_{\Gamma\Delta}^T]\left[\begin{array}{c}1_{p_I}\\
  1_{p_\Gamma}\end{array}\right].
\]
Therefore,  a basis of the one-dimensional null space of \EQ{bigeq} is
\begin{equation}
\label{equation:bignull}
\left( \begin{array}{cccccc}
0, & 1_{p_I}, & 0, & 0, & 1_{p_\Gamma}, &
-B_{\Delta,D}[B_{I\Delta}^T ~~ B_{\Gamma\Delta}^T]\left[\begin{array}{c}1_{p_I}\\ 1_{p_\Gamma}\end{array}\right] \end{array} \right),
\end{equation}
which is also valid if $p_\Gamma$ in \EQ{bigeq} is empty and then the block $1_{p_\Gamma}$ in \EQ{bignull} disappears.

For the case when the coarse level primal velocity space ${\bf W}_\Pi$ contains both the subdomain corner and edge-average variables such that the divergence free boundary
condition~\EQ{divergencefree} is satisfied,  we have $\int_{\Omega_i} \nabla \cdot {\bf w}_{\Delta}^{(i)}  = 0$, for all ${\bf w}_{\Delta}^{(i)} \in {\bf W}^{(i)}_{\Delta}$, which is in matrix form $\left(1_{p ^{(i)}_I}^TB_{I\Delta}^{(i)}+1_{p^{(i)}_\Gamma}^TB^{(i)}_{\Gamma\Delta}\right){\bf w}_{\Delta}^{(i)}=0$. As a result
\be
\label{equation:dfree}
[B^T_{I\Delta}~~B^T_{\Gamma \Delta}]\left[\begin{array}{c}1_{{p_I}}\\ 1_{{p_\Gamma}}\end{array}\right]=0.
\ee
Therefore,  when the divergence free boundary
condition~\EQ{divergencefree} is enforced and $p_\Gamma$ in
\EQ{bigeq} is empty, we know from
\EQ{dfree} and \EQ{bignull} that 
the leading four-by-four diagonal block in the coefficient
matrix of \EQ{bigeq} is singular and its null space consists of all
vectors with a constant pressure and zero velocity.

\section{Reduced system of linear equations}
\label{section:Gmatrix}

We first describe the FETI-DP formulation proposed in \cite{li12} and then use it as a framework to represent
the two previous FETI-DP algorithms studied in \cite{li05} and \cite{kim10}. Based on \EQ{bigeq}, denote
\be
\label{equation:AGtilde}
\widetilde{A} = \left[
\begin{array}{cccc}
A_{II}      & B_{II}^T       & A_{I \Delta}       & A_{I \Pi}       \\[0.8ex]
B_{II}      & 0              & B_{I \Delta}       & B_{I \Pi}       \\[0.8ex]
A_{\Delta I}& B_{I \Delta} ^T& A_{\Delta\Delta}   & A_{\Delta \Pi}  \\[0.8ex]
A_{\Pi I}   & B_{I \Pi}^T    & A_{\Pi \Delta}     & A_{\Pi \Pi}   \end{array}
\right], ~~
B_C=\left[
\begin{array}{cccc}
B_{\Gamma I} & 0 & B_{\Gamma \Delta} & B_{\Gamma \Pi} \\[0.8ex]
0            & 0 & B_{\Delta}        & 0              \end{array}
\right], ~~
f = \left[
\begin{array}{l}
{\bf f}_I        \\[0.8ex]
0                \\[0.8ex]
{\bf f}_{\Delta} \\[0.8ex]
{\bf f}_\Pi
\end{array}
\right].
\ee
The variables  $\left( {\bf u}_I, ~p_I, ~{\bf u}_{\Delta},~{\bf u}_{\Pi} \right)$ can be eliminated from \EQ{bigeq} and we obtain a Schur complement problem for the variables $\left(p_{\Gamma}, ~\lambda \right)$
\begin{equation}
\label{equation:spdG}
G \left[ \begin{array}{c}
p_\Gamma         \\[0.8ex]
\lambda
\end{array} \right] ~ = ~ g,
\end{equation}
where
\begin{equation}
\label{equation:Gmatrix}
G = B_C \widetilde{A}^{-1} B_C^T, \quad
g = B_C \widetilde{A}^{-1}  f.
\end{equation}

\begin{myremark}\label{remark:Asingular}
The only case that  $\widetilde{A}$ is singular in the algorithm is when $p_\Gamma$ in \EQ{bigeq} is empty and the divergence free boundary condition~\EQ{divergencefree} is enforced, as discussed at the end of Section \ref{section:DDM}. However the definitions in \EQ{Gmatrix} are still valid, since $f$ and columns of $B_C^T$ are  in the range of $\Atilde$ and the kernel of $\Atilde$ is a subspace of the kernel of $B_C$. For the simplicity of notation, we still use $\widetilde{A}^{-1}$ in \EQ{Gmatrix} and in other places of this paper to represent the solution of system of linear equations, not necessarily the inverse of $\widetilde{A}$.
\end{myremark}

We can see that $-G$ is the Schur complement of the coefficient matrix
of \EQ{bigeq} with respect to the last two row blocks:
\[
\left[ \begin{array}{cc} I & 0 \\[0.8ex] -B_C \widetilde{A}^{-1} & I \end{array} \right]
\left[ \begin{array}{cc} \widetilde{A} & B_C^T \\[0.8ex] B_C & 0 \end{array} \right]
\left[ \begin{array}{cc} I & - \widetilde{A}^{-1} B_C^T \\[0.8ex] 0  & I \end{array} \right] =
\left[ \begin{array}{cc} \widetilde{A} & 0 \\[0.8ex] 0 & -G \end{array} \right].
\]

If $\Atilde$ is nonsingular, then from the Sylvester law of inertia we can see that $G$ is symmetric positive semi-definite and its null space is derived from the null space of the original coefficient matrix of \EQ{bigeq}, cf.~\EQ{bignull}, and a basis is given by,
\[
\left( \begin{array}{cc} 1_{p_\Gamma}, & - B_{\Delta,D}[B_{I\Delta}^T ~~ B_{\Gamma\Delta}^T]\left[\begin{array}{c}1_{p_I}\\ 1_{p_\Gamma} \end{array} \right]
\end{array} \right).
\]
We denote $X = Q_\Gamma \bigoplus \Lambda$. The range space of $G$, denoted by $R_G$, is a subspace of $X$. $R_G$ is orthogonal to the null space of $G$ and thus has the form
\be
\label{equation:Grange}
R_G=\left\{  \left[ \begin{array}{c}
g_{p_\Gamma}         \\[0.8ex]
g_{\lambda}
\end{array} \right] \in X: g_{p_\Gamma}^T 1_{{p_\Gamma}} -
g_{\lambda}^T \left(B_{\Delta,D}[B_{I\Delta}^T ~~ B_{\Gamma\Delta}^T]\left[\begin{array}{c}1_{{p_I}}\\ 1_{{p_\Gamma}}\end{array}\right]\right)=0\right\}.
\ee

When $p_\Gamma$ in \EQ{bigeq} is empty and the divergence free
boundary condition~\EQ{divergencefree} is enforced, as discussed in
Remark \ref{remark:Asingular},
$\Atilde$ is singular and then $G$ becomes positive definite. The range space formula~\EQ{Grange} is still valid, cf. \EQ{dfree}, and in fact $R_G$ becomes the whole $\Lambda$.

In both cases, the restriction of $G$ to its range $R_G$ is positive definite.
The fact that the solution of \EQ{bigeq} always exists for any given $\left( {\bf f}_I, ~{\bf f}_{\Delta}, ~{\bf f}_{\Pi} \right)$ on the right-hand side implies that the solution of~\EQ{spdG} exits for any $g$ defined in \EQ{Gmatrix}. Therefore $g \in R_G$. When the conjugate gradient method is applied to solve \EQ{spdG} with zero initial guess, all the iterates are in the Krylov subspace generated by $G$ and $g$, which is also a subspace of $R_G$, and where the conjugate gradient method cannot break down. After obtaining $\left( p_{\Gamma}, ~\lambda \right)$ from solving \EQ{spdG}, the other components $\left( {\bf u}_I, ~p_I, ~{\bf u}_{\Delta}, ~{\bf u}_{\Pi} \right)$ in \EQ{bigeq} are obtained by back substitution.

The main computation of multiplying $G$ by a vector is the product of $\widetilde{A}^{-1}$ with a vector in the structure of $f$.  We denote
\[ A_{rr} = \left[ \begin{array}{ccc}
A_{II}       & B_{II}^T        & A_{I \Delta}     \\[0.8ex]
B_{II}       & 0               & B_{I \Delta}     \\[0.8ex]
A_{\Delta I} & B_{I \Delta} ^T & A_{\Delta\Delta} \end{array} \right] ,  \quad
A_{\Pi r} = A_{r \Pi}^T = \left[ A_{\Pi I}  \quad B_{I \Pi}^T  \quad A_{\Pi \Delta} \right], \quad  f_r = \left[ \begin{array}{l}
{\bf f}_I        \\[0.8ex]
0                \\[0.8ex]
{\bf f}_{\Delta} \end{array} \right],
\]
and define the Schur complement (coarse level problem)
\[
S_{\Pi} = A_{\Pi \Pi} - A_{\Pi r} A_{rr}^{-1} A_{r \Pi}.
\]
Then the product $\widetilde{A}^{-1}  f$
can be represented by
\be
\label{equation:im}
\left[ \begin{array}{c} A_{rr}^{-1} f_r \\[0.8ex] 0 \end{array} \right] ~ + ~
\left[ \begin{array}{c} -A_{rr}^{-1} A_{r \Pi} \\[0.8ex] I_\Pi \end{array}  \right] ~ S_{\Pi}^{-1} ~
\left({\bf f}_\Pi - A_{\Pi r} A_{rr}^{-1} f_r \right),
\ee
which requires solving two subdomain Neumann type problems and one coarse level problem.

The formulation of the reduced system \EQ{spdG} is valid for using
both discontinuous and continuous pressures in the algorithm, which is
discussed in the following based on whether $p_\Gamma$ is empty or
not. We note that when the discontinuous pressure is used in
the algorithm and $p_\Gamma$ is empty,  \EQ{spdG} is the same
system as those obtained in \cite{li05} and \cite{kim10}.

\subsection{$p_\Gamma$ is non-empty}
\label{subsection:cp}

This happens when a continuous pressure finite element space is used, where $p_\Gamma$ represents all the subdomain boundary pressure degrees of freedom shared by neighboring subdomains. Then~\EQ{spdG} is a system for both the subdomain boundary pressures and the Lagrange multipliers. $\Atilde$ in \EQ{Gmatrix} and $A_{rr}$ in \EQ{im} are both invertible and the coarse level problem operator $S_\Pi$ is symmetric positive definite.

$p_\Gamma$ can also be non-empty for the case of using discontinuous pressures. Since there are no pressure degrees of freedom shared by neighboring subdomains, it is free in the algorithm to choose any number of subdomain pressure variables as $p_\Gamma$. If $p_\Gamma$ contains at least one pressure degree of freedom from each subdomain, then $\Atilde$ in \EQ{Gmatrix} and  $A_{rr}$ in \EQ{im} are still invertible and $S_\Pi$ is symmetric positive definite.



\subsection{$p_\Gamma$ is empty}
\label{subsection:dp}

$p_\Gamma$ can be empty only when discontinuous pressures are used in the
finite element space since no
pressure degrees of freedom are shared by neighboring subdomains.

When $p_\Gamma$ is empty, \EQ{spdG} becomes a system for the
Lagrange multipliers only and it is the same equation as those
obtained in~\cite{li05} and \cite{kim10}.  The implementation of the product of
$\widetilde{A}^{-1}$ with a vector as specified in \EQ{im} is the
same as in \cite{kim10} and the resulting $S_\Pi$ is symmetric positive definite.  Kim {\it et. al.}
\cite{kim10} considered only the first choice of the coarse level primal velocity
space ${\bf W}_\Pi$, namely it contains only the subdomain corner
velocities, for which they proved in~\cite[Lemma 3.1]{kim10} that $\Atilde$ and $A_{rr}$ are
both invertible.  However, $\Atilde$ and $A_{rr}$ both become singular when the second choice of
$\vvec{W}_\Pi$ is used to enforce the divergence free boundary
condition~\EQ{divergencefree} (required for using the Dirichlet preconditioner, cf. Section \ref{section:Dirichlet}), even though their singularities do not affect the multiplication of $\widetilde{A}^{-1}$ by a vector; see Remark \ref{remark:Asingular}.

In \cite{li05}, the divergence free boundary condition~\EQ{divergencefree} is enforced for using the Dirichlet preconditioner and a different implementation of multiplying $\widetilde{A}^{-1}$ with a vector was used to avoid the singularity of $A_{rr}$ in \EQ{im}. There the subdomain constant pressures are pulled out from $Q_I$ to form another vector $p_0$ and $Q_I$ contains only subdomain interior pressures with zero average on each subdomain; $p_0$ is combined with the coarse level primal velocity variables to form the coarse level problem. More precisely, $\widetilde{A}$ in \EQ{AGtilde} is represented by
\[
\widetilde{A} = \left[
\begin{array}{ccccc}
A_{II}       & B_{II}^{-^T}       & A_{I \Delta}     & A_{I \Pi}      & B_{II}^{0^T}       \\[0.8ex]
B_{II}^-     & 0                  & B_{I \Delta}^-   & B_{I \Pi}^-    & 0                  \\[0.8ex]
A_{\Delta I} & B_{I \Delta}^{-^T} & A_{\Delta\Delta} & A_{\Delta \Pi} & B_{I \Delta}^{0^T} \\[0.8ex]
A_{\Pi I}    & B_{I \Pi}^{-^T}    & A_{\Pi \Delta}   & A_{\Pi \Pi}    & B_{I \Pi}^{0^T}    \\[0.8ex]
B_{II}^0     & 0                  & B_{I \Delta}^0   & B_{I \Pi}^0    & 0
\end{array}
\right],
\]
where, e.g., $B_{II}^-$ and $B_{II}^0$, represent blocks corresponding to the subdomain interior pressures with zero average and the subdomain constant pressures, respectively. Denote correspondingly
\[ A_{rr} = \left[ \begin{array}{ccc}
A_{II}       & B_{II}^{-^T}       & A_{I \Delta}     \\[0.8ex]
B_{II}^-     & 0                  & B_{I \Delta}^-   \\[0.8ex]
A_{\Delta I} & B_{I \Delta}^{-^T} & A_{\Delta\Delta} \end{array} \right] ,  \quad
A_{\Pi r} = A_{r \Pi}^T =
\left[ \begin{array}{ccc}
A_{\Pi I} & B_{I \Pi}^{-T} & A_{\Pi \Delta}\\[0.8ex]
B_{II}^0  & 0              & B_{I \Delta}^0
\end{array} \right], \quad  f_r = \left[ \begin{array}{l}
{\bf f}_I        \\[0.8ex]
0                \\[0.8ex]
{\bf f}_{\Delta} \end{array} \right],
\]
where $A_{rr}$ is invertible, and define the Schur complement (coarse level problem)
\[
S_{\Pi} = \left[ \begin{array}{cc} A_{\Pi \Pi} & B_{I \Pi}^{0^T}\\[0.8ex] B_{I \Pi}^0 & 0\end{array}\right] - A_{\Pi r} A_{rr}^{-1} A_{r \Pi},
\]
which is a saddle point problem. Then the product $\widetilde{A}^{-1} f$ can be represented by
\[
\left[ \begin{array}{c} A_{rr}^{-1} f_r \\[0.8ex] 0 \end{array} \right] ~ + ~
\left[ \begin{array}{c} -A_{rr}^{-1} A_{r \Pi} \\[0.8ex] I_\Pi \end{array}  \right] ~ S_{\Pi}^{-1} ~
\left(\left[\begin{array}{c} {\bf f}_\Pi \\[0.8ex] 0\end{array}\right]- A_{\Pi r} A_{rr}^{-1} f_r \right).
\]


\begin{myremark}\label{remark:method12} Even though the implementations of multiplying $\widetilde{A}^{-1}$ with a vector proposed in~\cite{li05} and \cite{kim10} are different, the same system \EQ{spdG} for the Lagrange multipliers is solved.
When equipped with the same type preconditioners, their convergence rates are the same.
\end{myremark}

\section{Some techniques}
\label{section:techniques}

We first define certain norms for several vector/function spaces. We denote
\be
\label{equation:Wtilde}
\vvec{\Wtilde} = {\bf W}_I \bigoplus \vvec{\Wtilde}_{\Gamma}.
\ee
For any $\vvec{w}$ in $\vvec{\Wtilde}$, denote its restriction to subdomain $\Omega_i$ by ${\bf w}^{(i)}$. A subdomain-wise $H^1$-seminorm can be defined for functions in $\vvec{\Wtilde}$ by
\[
|\vvec{w}|^2_{H^1} = \sum_{i=1}^N |\vvec{w}^{(i)}|^2_{H^1(\Omega_i)}.
\]

Several of the following lemmas have been proved in \cite{li12}; they are presented here for the completeness of this paper. We denote in \EQ{bigeq}
\begin{equation}\label{equation:Btilde}
\widetilde{B} = \left[ \begin{array}{ccc}
B_{II}      & B_{I \Delta}       & B_{I \Pi}      \\[0.8ex]
B_{\Gamma I}& B_{\Gamma \Delta}  & B_{\Gamma \Pi}
\end{array} \right].
\end{equation}
The following lemma on the stability of $\widetilde{B}$ is \cite[Lemma 5.1]{li12}.

\begin{mylemma}
\label{lemma:BtildeStability}
For any $\vvec{w} \in \vvec{\Wtilde}$ and $q \in Q$, $\left< \Btilde {\bf w}, q \right> \leq | \vvec{w} |_{H^1} \| q \|_{L^2}$.
\end{mylemma}

\beginproof \quad
$\displaystyle{\left<\Btilde {\bf w}, q \right>^2 = \left( \sum_{i=1}^N \int_{\Omega_i} \nabla\cdot {\bf w}^{(i)} q \right)^2\leq \left( \sum_{i=1}^N \sqrt{\int_{\Omega_i} | \nabla {\bf w}^{(i)} |^2} \sqrt{\int_{\Omega_i} q^2} \right)^2}$
\[
\qquad \le \left( \sum_{i=1}^N \int_{\Omega_i} | \nabla {\bf w}^{(i)} |^2 \right) \left( \sum_{i=1}^N \int_{\Omega_i} q^2 \right) = | \vvec{w} |^2_{H^1} \| q \|^2_{L^2}. \qquad \Box
\]

We define
\[
W = {\bf W}_I \bigoplus Q_I \bigoplus {\bf W}_{\Delta} \bigoplus {\bf W}_\Pi,
\]
and its subspace
\be
\label{equation:W0}
\Wtilde_{0} = \left\{ w = \left( {\bf w}_I, ~p_I, ~{\bf w}_{\Delta}, ~{\bf w}_{\Pi} \right) \in W :  B_{I I} \vvec{w}_I + B_{I\Delta} \vvec{w}_\Delta + B_{I\Pi} \vvec{w}_\Pi = 0 \right\}.
\ee
For any $w = \left( {\bf w}_I, ~p_I, ~{\bf w}_{\Delta}, ~{\bf w}_{\Pi} \right) \in \Wtilde_{0}$,
\be\label{equation:wg0}
\left< w, w \right>_{\widetilde{A}} = \sum_{i=1}^N \left[ \begin{array}{c} {\bf w}_I^{(i)} \\ {\bf w}_{\Delta}^{(i)} \\ {\bf w}_{\Pi}^{(i)} \end{array} \right]^T
\left[ \begin{array}{cccc}
A_{II}^{(i)}       & A_{I \Delta}^{(i)}      & A_{I \Pi}^{(i)}      \\[0.8ex]
A_{\Delta I}^{(i)} & A_{\Delta\Delta}^{(i)}  & A_{\Delta \Pi}^{(i)} \\[0.8ex]
A_{\Pi I}^{(i)}    & A_{\Pi \Delta}^{(i)}    & A_{\Pi \Pi}^{(i)}
\end{array} \right] \left[ \begin{array}{c} {\bf w}_I^{(i)} \\ {\bf w}_{\Delta}^{(i)} \\ {\bf w}_{\Pi}^{(i)} \end{array} \right]
= \sum_{i=1}^N \left| \left[ \begin{array}{c} {\bf w}_I^{(i)} \\ {\bf w}_{\Delta}^{(i)} \\ {\bf w}_{\Pi}^{(i)} \end{array} \right] \right|_{H^1(\Omega^i)}^2,
\ee
i.e., $\left< \cdot , \cdot \right>_{\Atilde}$ defines a seminorm on $\Wtilde_0$. In \EQ{wg0}, the superscript ${}^{(i)}$ is used to represent the restrictions of corresponding vectors and matrices to subdomain $\Omega_i$.

The following lemma is \cite[Lemma 6.6]{li12}.

\begin{mylemma}
\label{lemma:BcW0}
For any $w = \left( {\bf w}_I, ~p_I, ~{\bf w}_{\Delta}, ~{\bf w}_{\Pi} \right) \in \Wtilde_0$, $B_C w \in R_G$.
\end{mylemma}

\beginproof
We know for any $\left( {\bf f}_I, ~{\bf f}_{\Delta}, ~{\bf f}_{\Pi} \right) \in \vvec{W}_I \bigoplus \vvec{W}_\Delta \bigoplus \vvec{W}_\Pi$, $g$ defined in \EQ{Gmatrix} is in $R_G$.
For any $w = \left( {\bf w}_I, ~p_I, ~{\bf w}_{\Delta}, ~{\bf w}_{\Pi} \right) \in \Wtilde_0$, from the definition of $\Atilde$ in \EQ{AGtilde}, there always exists $\left( {\bf f}_I, ~{\bf f}_{\Delta}, ~{\bf f}_{\Pi} \right) \in \vvec{W}_I \bigoplus \vvec{W}_\Delta \bigoplus \vvec{W}_\Pi$, such that
\[
\Atilde w = \left[
\begin{array}{l}
{\bf f}_I        \\[0.8ex]
0                \\[0.8ex]
{\bf f}_{\Delta} \\[0.8ex]
{\bf f}_\Pi
\end{array}
\right], \quad \mbox{i.e.,} \quad w = \Atilde^{-1} \left[
\begin{array}{l}
{\bf f}_I        \\[0.8ex]
0                \\[0.8ex]
{\bf f}_{\Delta} \\[0.8ex]
{\bf f}_\Pi
\end{array}
\right].
\]
Taking such $\left( {\bf f}_I, ~{\bf f}_{\Delta}, ~{\bf f}_{\Pi} \right)$, $g$ defined in \EQ{Gmatrix} is $B_C w$.
$\qquad \Box$

From \EQ{AGtilde}, we denote the first row of $B_C$ by
\[
\widetilde{B}_{\Gamma} =
\left[ B_{\Gamma I} \quad 0  \quad B_{\Gamma \Delta}  \quad  B_{\Gamma \Pi} \right];
\]
for the second row, we denote the restriction from $W$ onto ${\bf W}_{\Delta}$ by $\widetilde{R}_{\Delta}$, such that
for any $w = \left( {\bf w}_I, ~p_I, ~{\bf w}_{\Delta}, ~{\bf w}_{\Pi} \right) \in W$, $\widetilde{R}_{\Delta} w = {\bf w}_{\Delta}$. Then $G$, defined in \EQ{Gmatrix}, can be written as the following two-by-two block structure
\begin{equation}
\label{equation:Gtwo}
G = \left[ \begin{array}{cc} G_{p_\Gamma p_\Gamma} & G_{p_\Gamma \lambda} \\[0.8ex] G_{\lambda p_\Gamma} & G_{\lambda \lambda} \end{array} \right],
\end{equation}
where
\[
G_{p_\Gamma p_\Gamma} = \widetilde{B}_{\Gamma} \widetilde{A}^{-1} \widetilde{B}_{\Gamma}^T, \qquad
G_{p_\Gamma \lambda} = \widetilde{B}_{\Gamma} \widetilde{A}^{-1} \widetilde{R}_{\Delta}^T B_{\Delta}^T,
\]
\[
G_{\lambda p_\Gamma} = B_{\Delta} \widetilde{R}_{\Delta} \widetilde{A}^{-1} \widetilde{B}_{\Gamma}^T, \qquad
G_{\lambda \lambda} = B_{\Delta} \widetilde{R}_{\Delta} \widetilde{A}^{-1} \widetilde{R}_{\Delta}^T B_{\Delta}^T.
\]

We  define a certain jump operator across the subdomain interface $\Gamma$. Let
$P_{D,1}: W \rightarrow W$, be defined by
\[
P_{D,1} = \widetilde{R}_{\Delta}^T B_{\Delta, D}^T B_{\Delta} \widetilde{R}_{\Delta},
\]
cf.~\cite{li06A}.
We can see that application of $P_{D,1}$ to a vector essentially
computes the difference (jump) of the dual velocity components across
the subdomain interface and then distributes the jump to neighboring
subdomains according to the scaling factor $\delta^\dagger(x)$. In
fact, from the definition of $P_{D,1}$, the only component involved in
its application is the component in the space $\vvec{W}_\Delta$; all
other components are kept zero and they are added into the definition
to make $P_{D,1}$ more convenient to use in the analysis.

Note that for any $w \in W$, $\left< P_{D,1} w, P_{D,1} w \right>_{\widetilde{A}} = \left< B_{\Delta, D}^T B_{\Delta} {\bf w}_{\Delta}, B_{\Delta, D}^T B_{\Delta} {\bf w}_{\Delta} \right>_{A_{\Delta \Delta}}$.
The following lemma can be found essentially from \cite[Section 6]{li07}.
\begin{mylemma}
\label{lemma:jump}
There exists a function $\Phi_1(H/h)$, such that $\left< P_{D,1} w, P_{D,1} w \right>_{\widetilde{A}} \leq  \Phi_1(H/h) \left< w, w \right>_{\widetilde{A}}$,  for all $w \in \Wtilde_0$. For two-dimensional problems, when the coarse level primal velocity space contains the subdomain corner velocity variables,
$\Phi_1(H/h) = C (H/h) (1 + \log{(H/h)})$.
\end{mylemma}

To improve the bound on the jump operator, the jumps across the subdomain interface can be extended
to the interior of subdomains by subdomain discrete harmonic extension. We define a Schur complement $H^{(i)}_{\Delta}: {\bf W}^{(i)}_{\Delta} \rightarrow {\bf W}^{(i)}_{\Delta}$ by, for any ${\bf w}_{\Delta}^{(i)} \in {\bf W}^{(i)}_{\Delta}$,
\begin{equation}
\label{equation:SiDelta}
\left[ \begin{array}{cc}
A_{II}^{(i)}            & A_{I \Delta}^{(i)}     \\[0.8ex]
A_{\Delta I}^{(i)} & A_{\Delta\Delta}^{(i)}
\end{array} \right] \left[ \begin{array}{c}
{\bf w}_I^{(i)}        \\[0.8ex]
{\bf w}_{\Delta}^{(i)}
\end{array} \right] =
\left[ \begin{array}{l}
{\bf 0}          \\[0.8ex]
H_{\Delta}^{(i)} {\bf w}_\Delta^{(i)}
\end{array} \right] \mbox{ . }
\end{equation}
We can see that to multiply $H_{\Delta}^{(i)}$ by the vector $\vvec{w}_\Delta^{(i)}$, subdomain elliptic problems with given boundary velocity $\vvec{w}_\Delta^{(i)}$ and $\vvec{w}^{(i)}_\Pi=\vvec{0}$ need be solved. Using $H^{(i)}_{\Delta}$, we define the second jump operator $P_{D,2}: W \rightarrow W$, by: for any given $w = \left( {\bf w}_I, ~p_I, ~{\bf w}_{\Delta}, ~{\bf w}_{\Pi} \right) \in W$, on each subdomain $\Omega_i$, the subdomain interior velocity part of $P_{D,2} w$ is taken as ${\bf w}_I^{(i)}$ in the solution of \EQ{SiDelta}, with given subdomain boundary velocity ${\bf w}_{\Delta}^{(i)} = B_{\Delta, D}^{(i)^T} B_{\Delta} {\bf w}_{\Delta}$. Here $B_{\Delta, D}^{(i)^T}$ represents restriction of $B_{\Delta, D}^T$ on subdomain $\Omega_i$ and is a map from $\Lambda$ to $\vvec{W}^{(i)}_\Delta$. The other components of $P_{D,2} w$ are kept zero. Therefore
\begin{eqnarray}
\left<P_{D,2} w,P_{D,2} w\right>_{\Atilde} & = & \sum_{i=1}^N
\left[ \begin{array}{c}
{\bf w}_I^{(i)}        \\[0.8ex]
{\bf w}_{\Delta}^{(i)}
\end{array} \right]^T  \left[ \begin{array}{ccc}
A_{II}^{(i)}          & A_{I \Delta}^{(i)}     \\[0.8ex]
A_{\Delta I}^{(i)}& A_{\Delta\Delta}^{(i)}
\end{array} \right] \left[ \begin{array}{c}
{\bf w}_I^{(i)}        \\[0.8ex]
{\bf w}_{\Delta}^{(i)}
\end{array} \right]  \nonumber \\
& = &  \sum_{i=1}^N  {\bf w}_\Delta^{(i)^T} H_{\Delta}^{(i)} {\bf w}_\Delta^{(i)} = \sum_{i=1}^N {\bf w}_{\Delta}^T B_{\Delta}^T
B_{\Delta, D}^{(i)} H_{\Delta}^{(i)} B_{\Delta, D}^{(i)^T} B_{\Delta}
{\bf w}_{\Delta} \label{equation:PD2}  \\
& = & \sum_{i=1}^N \left|\left[\begin{array}{c} B_{\Delta,
        D}^{(i)^T} B_{\Delta} {\bf w}_{\Delta}\\
      0 \end{array}\right]\right|_{H^{1/2}(\partial \Omega^i)}^2 \leq
\Phi_2(H/h) \sum_{i=1}^N |\vvec{w}^{(i)}_\Gamma|^2_{H^{1/2}(\partial \Omega^i)},  \nonumber
\end{eqnarray}
where $\vvec{w}^{(i)}_\Gamma$ represents the restriction of $\left( {\bf w}_{\Delta}, ~{\bf w}_{\Pi} \right)$ on subdomain $\Omega_i$. The last inequality in \EQ{PD2} is a well established result, cf., \cite[Lemma 6.34]{Toselli:2004:DDM}. It has been established that, $\Phi_2(H/h) = C (1 + \log{(H/h)})^2$, for two-dimensional problems, when the coarse level primal velocity space contains the subdomain corner velocity variables, cf.~\cite[Lemma 3.3]{4wid}.

Then from \EQ{PD2} and \EQ{wg0}, we have
\begin{mylemma}
\label{lemma:jumptwo}
There exists a function $\Phi_2(H/h)$, such that $\left< P_{D,2} w, P_{D,2} w \right>_{\widetilde{A}}\leq \Phi_2(H/h) \left< w, w \right>_{\widetilde{A}}$, for all $w \in \Wtilde_0$. For two-dimensional problems, when the coarse level primal velocity space contains the subdomain corner velocity variables, $\Phi_2(H/h) = C (1 + \log{(H/h)})^2$.
\end{mylemma}

The following lemma is also used in our analysis and can be found at \cite[Lemma~2.3]{paulo2003}.
\begin{mylemma}
\label{lemma:paul}
Consider the saddle point problem: find  $(\vvec{u}, p) \in \vvec{W} \bigoplus Q$, such that
\be \left[
\begin{array}{cc}
A       & B^T      \\[0.8ex]
B       & 0
\end{array}
\right] \left[
\begin{array}{l}
\vvec{u}        \\[0.8ex]
p
\end{array}
\right]
=\left[
\begin{array}{l}
\vvec{f}        \\[0.8ex]
g
\end{array}
\right],
\ee
where $A$ and $B$ are as in \EQ{matrix}, $\vvec{f} \in \vvec{W}$, $g \in Im(B) \subset Q$. Let $\beta$ be the inf-sup constant specified in \EQ{infsup}.
Then $\| \vvec{u} \|_A \le \| \vvec{f} \|_{A^{-1}} + \frac{1}{\beta} \| g \|_{Z^{-1}}$,
where $Z$ is the mass matrix defined in Section~\ref{section:FEM}.
\end{mylemma}

The lumped and the Dirichlet preconditioners for solving \EQ{spdG} will be studied in Sections~\ref{section:lumped} and \ref{section:Dirichlet}. Let $M^{-1}$ be symmetric positive definite. When the conjugate gradient method is applied to solve the preconditioned system
\begin{equation}
\label{equation:Mspd}
M^{-1} G x ~ = ~ M^{-1} g,
\end{equation}
with zero initial guess, all iterates belong to the Krylov subspace generated by the operator $M^{-1} G$ and the vector $M^{-1}
g$, which is a subspace of the range of $M^{-1} G$. We denote the
range of $M^{-1} G$ by $R_{M^{-1} G}$. The following two lemmas are also from \cite{li12}.

\begin{mylemma}
\label{lemma:CG} Let $M^{-1}$ be symmetric positive definite.
Then the conjugate gradient method applied to solving \EQ{Mspd} with zero initial guess cannot break down.
\end{mylemma}

\beginproof
We just need to show that for any $0 \neq x \in R_{M^{-1} G}$, $G x \neq 0$. Let $0 \neq x = M^{-1} G y$, for a certain $y \in X$ and $y \neq 0$. Then $G x = G M^{-1} G y$, which cannot be zero since $G y \neq 0$ and $y^T  G M^{-1} G y \neq 0$.
\endproof

\begin{mylemma}
\label{lemma:m1R} Let $M^{-1}$ be symmetric positive definite. For any $x \in R_{M^{-1} G}$,
\[
\left<Mx, x \right> = \max_{y \in R_G, y \neq 0} \frac{\left<y, x\right>^2}{\left<M^{-1}y,y\right>}.
\]
\end{mylemma}
\beginproof Denote the range of $M^{-\frac12} G$ by $R_{M^{-1/2} G}$,
then for any $\lambda \in R_{M^{-1} G}$
\begin{eqnarray*}
\left<Mx, x \right> & = & \left<M^{\frac12} x, M^{\frac12} x \right> =
\max_{z \in R_{M^{-1/2} G}, z \neq 0} \frac{\left<M^{\frac12} x, z\right>^2}{\left<z, z\right>} \\
& = & \max_{y \in R_G, y \neq 0}
\frac{\left<M^{\frac12} x, M^{-\frac12} y\right>^2}{\left<M^{-\frac12} y, M^{-\frac12} y\right>}
= \max_{y \in R_G, y \neq 0} \frac{\left<y,x\right>^2}{\left<M^{-1}y,y\right>}. \qquad \Box
\end{eqnarray*}

\section{The lumped preconditioner}
\label{section:lumped}

The lumped  preconditioner for solving \EQ{spdG} has been studied in \cite{li12} for the case of using continuous pressures in the finite element discretization. The preconditioner is given by
\[
M_{L}^{-1} = \left[ \begin{array}{cc} \frac{1}{h^2} I_{p_\Gamma} & \\[0.8ex]
& M^{-1}_{L,\lambda}  \end{array} \right],
\]
where $I_{p_\Gamma}$ is the identity matrix of the same length as $p_\Gamma$,
\[
M^{-1}_{L, \lambda} = B_{\Delta, D} \widetilde{R}_{\Delta} \widetilde{A} \widetilde{R}_{\Delta}^T B_{\Delta, D}^T.
\]

It has been proved in \cite{li12} that,  for the case of using
continuous pressures, the condition number of the preconditioned
operator $M_{L}^{-1} G$ be bounded by $C \Phi_1(H,h)/\beta^2$, where
$\Phi_1(H,h)$ is defined in Lemma \LA{jump} and $\beta$ is the inf-sup
constant specified in~\EQ{infsup}. The same bound also holds for more
general cases where $p_\Gamma$ is non-empty, including for using discontinuous pressures; see Section~\ref{subsection:cp} and~\cite[Remark~6.10]{li12}.

When $p_\Gamma$ is empty, \EQ{spdG} becomes a system for the Lagrange
multipliers only, which is the same as the
one obtained in \cite{kim10}, cf. Section \ref{subsection:dp}. The matrix $G$ in \EQ{spdG} contains only the second diagonal block in \EQ{Gtwo}, i.e., \[
G = G_{\lambda \lambda} = B_{\Delta} \widetilde{R}_{\Delta} \widetilde{A}^{-1} \widetilde{R}_{\Delta}^T B_{\Delta}^T
\]
and the lumped preconditioner $M^{-1}_{L}$ simplifies to
\[
M^{-1}_{L} = M^{-1}_{L, \lambda} = B_{\Delta, D} \widetilde{R}_{\Delta} \widetilde{A} \widetilde{R}_{\Delta}^T B_{\Delta, D}^T.
\]

In this section, we establish the condition number bound for the case of empty
$p_\Gamma$ following the approach given in \cite{li12}.  The bound we
obtained here is the same as that obtained in \cite{kim10}, but this approach  greatly simplifies the analysis given in \cite{kim10}. We also note that this bound is also the same as that obtained for non-empty $p_\Gamma$
in \cite{li12}.

As discussed in \cite{kim10} and \cite{li12}, for using the lumped
preconditioner, we only need use the first choice of the coarse
level primal velocity space, i.e., ${\bf W}_{\Pi}$ contains only the subdomain corner velocities. We also note that when $p_\Gamma$ is empty, $p_I$ contains all the pressure degrees of freedom and the results in Section \ref{section:techniques} are still valid. In fact the analysis given in the following is essentially the special case of the analysis in \cite{li12} when the blocks corresponding to $p_\Gamma$ no longer exist. We have the following lemmas.

\begin{mylemma}\label{lemma:upperFlump}
For any $w\in \Wtilde_{0}$, $\left<M^{-1}_{L} B_C w, B_C w \right>\le \Phi_1(H,h)\left<\Atilde w, w\right>$,
where $\Phi_1(H,h)$ is as defined in Lemma \LA{jump}.
\end{mylemma}

\beginproof
For any given $w \in \Wtilde_{0}$, we have
\BA
\label{equation:MBWF}
\left<M^{-1}_{L} B_C w, B_C w\right> & = & \left(B_\Delta\Rtilde_\Delta w\right)^T M^{-1}_{L}B_\Delta\Rtilde_\Delta w\nonumber\\
&=&\left(B_\Delta\Rtilde_\Delta w\right)^T
B_{\Delta,D}\Rtilde_{\Delta}\Atilde\Rtilde_{\Delta}^TB_{\Delta,D}^T
\left(B_\Delta\Rtilde_\Delta w\right)\nonumber\\
&=&\left<P_{D,1}w,P_{D,1}w\right>_{\Atilde} \le\Phi_1(H,h)\left<w,w\right>_{\Atilde},
\EA
where we have used Lemma \LA{jump} for the last inequality.
\endproof

\begin{mylemma}
\label{lemma:lowerFlump}
For any given $y =  g_\lambda \in R_G$, there exits
$w \in \Wtilde_{0}$, such that $B_C w = y$, and $\left <\Atilde w, w\right> \le \frac{C}{\beta^2}  \left< M_{L}^{-1}y, y \right>$.
\end{mylemma}

\beginproof
Given $y =g_\lambda \in R_G$, let ${\bf u}_{\Delta}^{(I)} = B_{\Delta, D}^T g_\lambda$, $\u^{(I)} = \left({\bf 0}, ~{\bf u}_{\Delta}^{(I)}, {\bf 0}\right) \in {\bf W}_I \bigoplus {\bf W}_{\Delta} \bigoplus {\bf W}_\Pi$, and $w^{(I)} = \left( {\bf 0}, ~0, ~{\bf u}_{\Delta}^{(I)}, ~ {\bf 0} \right) \in {\bf W}_I \bigoplus Q_I \bigoplus {\bf W}_{\Delta} \bigoplus {\bf W}_\Pi$.  We have
\be
\label{equation:uOne}
| \u^{(I)} |^2_{H^1} =  \left< A_{\Delta\Delta}\u^{(I)}_\Delta,\u^{(I)}_\Delta\right>,
\ee
and
\be
\label{equation:bcWone}
B_C w^{(I)} = \left[ \begin{array}{cccc}
0            & 0 & B_{\Delta}        & 0              \end{array}
\right] \left[ \begin{array}{c}
{\bf 0} \\[0.8ex] 0 \\[0.8ex] B_{\Delta, D}^T g_\lambda \\[0.8ex] {\bf 0}
\end{array} \right] = g_\lambda ,
\ee
where we used the fact that $B_\Delta B_{\Delta, D}^T = I$.

We consider a solution to the following fully assembled system of linear equations of the form~\EQ{matrix}: find $\left( {\bf u}_I^{(II)}, ~p_I^{(II)}, ~{\bf u}_{\Gamma}^{(II)} \right) \in {\bf W}_I \bigoplus Q_I \bigoplus {\bf W}_{\Gamma}$, such that
\be
\label{equation:uTwo}
\left[
\begin{array}{ccc}
A_{II}      & B_{II}^T       & A_{I \Gamma}       \\[0.8ex]
B_{II}      & 0              & B_{I \Gamma}                      \\[0.8ex]
A_{\Gamma I}& B_{I \Gamma} ^T& A_{\Gamma\Gamma}  \\[0.8ex]
\end{array}
\right]
\left[ \begin{array}{c}
{\bf u}_I^{(II)}        \\[0.8ex]
p_I^{(II)}              \\[0.8ex]
{\bf u}_{\Gamma}^{(II)} \\[0.8ex]
\end{array} \right] =
\left[ \begin{array}{l}
{\bf 0}        \\[0.8ex]
-B_{I\Delta}{\bf u}^{(I)}_\Delta              \\[0.8ex]
{\bf 0}        \\[0.8ex]
\end{array} \right] \mbox{ ,  }
\ee
where a particular right-hand side is chosen.
We first note that, since $g_\lambda \in R_G$, the right-hand side vector of the above system satisfies, cf.~\EQ{Grange},
\[
1_{p_I}^T \left( -B_{I\Delta}{\bf u}^{(I)}_\Delta \right)
=  - 1_{p_I}^T B_{I\Delta} B_{\Delta, D}^T g_\lambda
= 0,
\]
i.e., it has a zero average, which implies the existence of solutions to the above system. Recall that here 
$p_I$ contains all pressure degrees of freedom and there are no subdomain boundary pressure degrees of freedom.

Denote ${\bf u}^{(II)} = \left( {\bf u}_I^{(II)}, ~{\bf u}_{\Gamma}^{(II)} \right) \in {\bf W}$. From the inf-sup stability of the original problem~\EQ{matrix} and Lemma \ref{lemma:paul}, we have
\be \label{equation:uIboundF}
| {\bf u}^{(II)} |^2_{H^1} \leq
\frac{1}{\beta^2} \left\| B_{I\Delta}\u^{(I)}_\Delta\right\|^2_{Z^{-1}}.
\ee

The right-hand side of \EQ{uIboundF} can be bounded, using Lemma \LA{BtildeStability}, as follows,
\begin{eqnarray}\label{equation:uuOneF}
&&\left\| B_{I\Delta}\u^{(I)}_\Delta \right\|^2_{Z^{-1}}=\left<\Btilde {\bf u}^{(I)},  \Btilde {\bf u}^{(I)} \right>_{Z^{-1}} =  C \max_{q \in Q} \frac{\left<\Btilde {\bf u}^{(I)}, q \right>^2}{\left<q, q\right>_Z} \nonumber\\
& \le & C \max_{q \in Q} \frac{| {\bf u}^{(I)} |_{H^1}^2  \| q \|^2_{L^2}}{\| q \|^2_{L^2}} = C \left< A_{\Delta\Delta}\u^{(I)}_\Delta,\u^{(I)}_\Delta\right>.
\end{eqnarray}

Split the continuous subdomain boundary velocity ${\bf u}_{\Gamma}^{(II)}$ into the dual part ${\bf u}_{\Delta}^{(II)} \in {\bf W}_{\Delta}$ and the primal part ${\bf u}_{\Pi}^{(II)} \in {\bf W}_{\Pi}$, and denote $w^{(II)} = \left( {\bf u}_I^{(II)}, ~p_I^{(II)}, ~{\bf u}_{\Delta}^{(II)}, ~{\bf u}_{\Pi}^{(II)} \right)$.
We have, from~\EQ{uTwo},
\be
\label{equation:biWtwo}
\left[ \begin{array}{cccc} B_{II} & 0 & B_{I \Delta} & B_{I \Pi} \end{array} \right]
\left[ \begin{array}{c}
{\bf u}_I^{(II)}        \\[0.8ex]
p_I^{(II)}              \\[0.8ex]
{\bf u}_{\Delta}^{(II)} \\[0.8ex]
{\bf u}_{\Pi}^{(II)}
\end{array} \right] = -B_{I\Delta}\u^{(I)}_\Delta,
\ee
and
\be
\label{equation:bcWtwo}
B_C w^{(II)} = \left[
\begin{array}{cccc}
0            & 0 & B_{\Delta}        & 0              \end{array}
\right] \left[ \begin{array}{c}
{\bf u}_I^{(II)}        \\[0.8ex]
p_I^{(II)}              \\[0.8ex]
{\bf u}_{\Delta}^{(II)} \\[0.8ex]
{\bf u}_{\Pi}^{(II)}
\end{array} \right] = 0.
\ee

Let $w = w^{(I)} + w^{(II)}$. We can then see from \EQ{biWtwo} that $w \in \Wtilde_{0}$, cf. \EQ{W0}. We can also see from \EQ{bcWone} and \EQ{bcWtwo} that $B_C w = y$. Furthermore, by \EQ{wg0},
\[
| w |^2_{\widetilde{A}}=| \u^{(I)}+ \u^{(II)}|^2_{H^1}  \le| \u^{(I)} |^2_{H^1}+| \u^{(II)} |^2_{H^1} \leq
\frac{C}{\beta^2} \left< A_{\Delta\Delta}\u^{(I)}_\Delta,\u^{(I)}_\Delta\right>,
\]
where we have used \EQ{uOne}, \EQ{uIboundF} and \EQ{uuOneF} for the last inequality.

On the other hand, we have
\begin{eqnarray*}
\left< M_{L}^{-1}y,y \right> & = &
g_{\lambda}^T M^{-1}_{L}g_\lambda
= g_{\lambda}^T
B_{\Delta,D}\Rtilde_{\Delta}\Atilde\Rtilde_{\Delta}^TB_{\Delta,D}^Tg_{\lambda}
= \left<A_{\Delta\Delta} \vvec{u}^{(I)}_\Delta, \vvec{u}^{(I)}_\Delta\right>. \qquad \Box
\end{eqnarray*}

\begin{mytheorem}
\label{theorem:tcondF} For all $ x=~\lambda\in R_{M_{L}^{-1} G}$,
\[
c \beta^2 \left<M_{L}x,x \right>\leq \left< G x,x \right> \leq \Phi_1(H,h) \left<
M_{L}x,x \right>,
\]
where $\Phi_1(H,h)$ is as defined in Lemma \LA{jump}, $\beta$ is the inf-sup constant specified in \EQ{infsup}.
\end{mytheorem}

\beginproof
$\quad \displaystyle{\left< Gx,x\right> = x^T B_C \Atilde^{-1} B_C^Tx= x^T B_C \Atilde^{-1} \Atilde \Atilde^{-1} B_C^Tx =
\left< \Atilde^{-1} B_C^Tx, \Atilde^{-1} B_C^Tx\right>_{\Atilde}}$.

\noindent Since $\Atilde^{-1} B_C^Tx \in \Wtilde_{0}$ and $\left< \cdot, \cdot \right>_{\Atilde}$ defines an inner product on $\Wtilde_{0}$, we have
\be \label{equation:Fnorm}
\left< Gx,x\right> = \max_{v \in \Wtilde_{0}, v \neq 0} \frac{\left<v, \Atilde^{-1} B_C^Tx \right>_{\widetilde{A}}}{\left<v, v \right>_{\widetilde{A}}} =\max_{v \in \Wtilde_0, v \neq 0}
\frac{\left<B_Cv,x\right>^2}{\left<\Atilde v,v\right>}.
\ee

{\it Lower bound:} From Lemma \ref{lemma:lowerFlump}, we know that for any given $y = g_\lambda \in R_G$, there exits
$w \in \Wtilde_{0}$, such that $B_C w = y$ and $\left <\Atilde w, w\right> \le \frac{C}{\beta^2}  \left< M_{L}^{-1}y, y \right>$. Then from  \EQ{Fnorm}, we have
$$
\left< Gx,x\right> \ge\frac{\left<B_C w,x\right>^2}{\left<\Atilde w,w\right>}
\ge c \beta^2 \frac{\left<y,x\right>^2}{\left<M_{L}^{-1}y,y\right>}.
$$
Since $y$ is arbitrary, using Lemma \ref{lemma:m1R}, we have
$$\left< Gx,x\right> \ge c  \beta^2 \max_{y \in R_G, y \neq 0} \frac{\left<y,x\right>^2}{\left<M_{L}^{-1}y,y\right>} = c \beta^2 \left<M_{L}x,x \right>.$$

{\it Upper bound:} From \EQ{Fnorm}, Lemmas \ref{lemma:BcW0}, \ref{lemma:upperFlump} and \ref{lemma:m1R},  we have
\begin{eqnarray*}
\left< Gx,x\right> & \le & \Phi_1(H,h)\max_{v\in \Wtilde_{0}, v \neq 0} \frac{\left<B_C v,x\right>^2}{\left<M^{-1}_{L} B_C v,B_C v\right>} \\
& \le & \Phi_1(H,h) \max_{y\in R_G, y \neq 0} \frac{\left<y,x\right>^2}{\left<M^{-1}_{L} y,y\right>} = \Phi_1(H,h) \left<M_{L}x,x\right>. \qquad \Box
\end{eqnarray*}

\section{The Dirichlet  preconditioner}
\label{section:Dirichlet}

In the lumped preconditioner discussed in the previous section, the subdomain interface jump of the velocity component is extended by zero to the interior of subdomains. Comparing Lemmas \LA{jump} and \LA{jumptwo}, the discrete subdomain harmonic extension of the jump to the interior of subdomains has a better stability, which leads to the Dirichlet preconditioner discussed in this section.

\begin{myremark}
Subdomain discrete Stokes extensions, obtained by solving saddle-point problems,
are used for the Dirichlet
preconditioner studied in \cite{li05} for solving \EQ{spdG} with
discontinuous pressures. In the Dirichlet preconditioner proposed in
this section, the discrete subdomain harmonic extensions of the jump,
obtained by solving symmetric positive definite problems,  are used in each
iteration. Even though these two extensions are spectrally equivalent,
using the subdomain harmonic extensions is more efficient  than solving indefinite subdomain problems.
\end{myremark}

We define a Schur complement $S^{(i)}_{\Delta\Pi}: {\bf W}^{(i)}_{\Delta} \bigoplus  {\bf W}^{(i)}_{\Pi} \rightarrow {\bf W}^{(i)}_{\Delta} \bigoplus  {\bf W}^{(i)}_{\Pi}$. For any ${\bf u}_{\Gamma}^{(i)} = \left( {\bf u}_{\Delta}^{(i)}, ~ {\bf u}_{\Pi}^{(i)} \right) \in {\bf W}^{(i)}_{\Delta} \bigoplus  {\bf W}^{(i)}_{\Pi}$, $S^{(i)}_{\Delta\Pi} {\bf u}_{\Gamma}^{(i)}$ is determined by
\be \label{equation:schurplus}
\left[
\begin{array}{ccccc}
A_{II}^{(i)}      & B_{II}^{(i)^T}      & A_{I \Delta}^{(i)}      &A_{I \Pi}^{(i)}      &B_{\Gamma I}^{(i)^T} \\[0.8ex]
B_{II}^{(i)}      & 0                   & B_{I \Delta}^{(i)}      &B_{I \Pi}^{(i)}
&0                    \\[0.8ex]
A_{\Delta I}^{(i)}& B_{I \Delta}^{(i)^T}& A_{\Delta\Delta}^{(i)}  &A_{\Delta \Pi}^{(i)} &B_{\Gamma \Delta}^{(i)^T} \\[0.8ex]
A_{\Pi I}^{(i)}   & B_{I \Pi}^{(i)^T}   & A_{\Pi \Delta}^{(i)}    &A_{\Pi \Pi}^{(i)}    &B_{\Gamma \Pi}^{(i)^T}    \\[0.8ex]
B_{\Gamma I}^{(i)}& 0                   & B_{\Gamma \Delta}^{(i)} &B_{\Gamma \Pi}^{(i)}
&0
\end{array}
\right]
\left[ \begin{array}{c}
{\bf u}_I^{(i)}        \\[0.8ex]
p_I^{(i)}              \\[0.8ex]
{\bf u}_{\Delta}^{(i)} \\[0.8ex]
{\bf u}_{\Pi}^{(i)}    \\[0.8ex]
p_\Gamma^{(i)}
\end{array} \right] =
\left[ \begin{array}{l}
{\bf 0}          \\[0.8ex]
0                \\[0.8ex]
\left( S^{(i)}_{\Delta\Pi} {\bf u}_{\Gamma}^{(i)} \right)_\Delta  \\[0.8ex]
\left( S^{(i)}_{\Delta\Pi} {\bf u}_{\Gamma}^{(i)} \right)_\Pi   \\[0.8ex]
0
\end{array} \right] \mbox{ . }
\ee

The solution of \EQ{schurplus} will be needed below in the proof of
Lemma \LA{lowerGDiri}. The subdomain incompressible Stokes problem
\EQ{schurplus} contains the subdomain constant pressure component. To
guarantee the existence of its solution, the given subdomain boundary
velocity ${\bf u}_{\Gamma}^{(i)} = \left( {\bf u}_{\Delta}^{(i)}, ~
  {\bf u}_{\Pi}^{(i)} \right)$ need satisfy the divergence free
condition $\int_{\partial \Omega_i} {\bf u}_{\Gamma}^{(i)} \cdot {\bf
  n} = 0$. In our applications, ${\bf u}_{\Gamma}^{(i)}$ will
represent the jump of the subdomain interface velocity across the
subdomain boundary and the coarse level component ${\bf
  u}_{\Pi}^{(i)}$ will always be zero. The divergence free
condition will be satisfied if \EQ{divergencefree} is enforced. Due
to this reason, for the condition number bound analysis for using the Dirichlet preconditioner, we assume that the
second choice of the coarse level primal velocity space ${\bf W}_\Pi$ as described in Section \ref{section:DDM}
is used.

\begin{myremark}
\label{remark:constraints}
The assumption that \EQ{divergencefree} is enforced is only required for the condition number bound analysis. When \EQ{divergencefree} is not enforced, e.g., when ${\bf W}_\Pi$ contains only the subdomain corner velocity variables, the Dirichlet preconditioner discussed below is still symmetric positive definite and the preconditioned conjugate gradient method can still be used, even though its condition number bound is no longer available.
\end{myremark}

The following lemma is a well established result, cf. \cite[Theorem
4.1]{3brpa} or \cite[Lemma 3.1]{8pavwid}.
\begin{mylemma}
\label{lemma:stokesnormPlus} For all ${\bf u}_{\Gamma}^{(i)} = \left( {\bf u}_{\Delta}^{(i)}, ~ {\bf u}_{\Pi}^{(i)} \right) \in {\bf W}^{(i)}_{\Delta} \bigoplus  {\bf W}^{(i)}_{\Pi}$,
$$
c \beta | {\bf u}_{\Gamma}^{(i)} |_{S_{\Delta\Pi}^{(i)}} \leq | {\bf u}_{\Gamma}^{(i)} |_{H^{1/2}(\partial \Omega^i)} \leq | {\bf u}_{\Gamma}^{(i)} |_{S_{\Delta\Pi}^{(i)}} ,
$$
where $\beta$ is the inf-sup constant specified in \EQ{infsup}.
\end{mylemma}

We also denote the direct sum of the discrete subdomain harmonic extension operators $H^{(i)}_\Delta, i = 1, \ldots, N,$ defined in \EQ{SiDelta}, by $H_{\Delta}: {\bf W}_{\Delta} \rightarrow {\bf W}_{\Delta}$.

In the following, the condition number bound for using the Dirichlet preconditioner is established for the case when $p_\Gamma$ is non-empty. As discussed in Section \ref{subsection:cp}, $p_\Gamma$ in \EQ{spdG} is non-empty when either the continuous pressure is used in the finite element discretization, or the discontinuous pressure is used and $p_\Gamma$ contains at least one pressure degree of freedom from each subdomain. In fact, the same condition number bound also holds for the case when $p_\Gamma$ is empty, which will be discussed briefly at the end of this section.

We define $M^{-1}_{D,\lambda}$ by
\be
\label{equation:subdomainwise}
M^{-1}_{D,\lambda} =B_{\Delta, D} H_{\Delta} B_{\Delta, D}^.
\ee
and the Dirichlet preconditioner for solving \EQ{spdG} is
\[
M_{D}^{-1} = \left[ \begin{array}{cc} \frac{1}{h^2} I_{p_\Gamma} & \\[0.8ex]
& M^{-1}_{D,\lambda}  \end{array} \right].
\]

\begin{mylemma}\label{lemma:upperGDiri}
For any $w\in \Wtilde_{0}$, $\left<M^{-1}_{D} B_C w,B_Cw\right>\le C\Phi_2(H,h)\left<\Atilde w, w\right>$,
where $\Phi_2(H,h)$ is as defined in Lemma~\LA{jumptwo}.
\end{mylemma}

\beginproof
Given $w = \left( {\bf w}_I, ~p_I, ~{\bf w}_{\Delta}, ~{\bf w}_{\Pi} \right) \in \Wtilde_{0}$, let
$g_{p_\Gamma} = B_{\Gamma I} {\bf w}_I + B_{\Gamma \Delta} {\bf w}_\Delta + B_{\Gamma\Pi} {\bf w}_\Pi$. Similar to Lemma
\LA{upperFlump}, we have from \EQ{PD2},
\begin{eqnarray*}
& & \left<M^{-1}_{D} B_C w,B_C w\right> =
\frac{1}{h^2}\left<g_{p_\Gamma}, g_{p_\Gamma} \right>+ \left(B_\Delta {\bf w}_\Delta
\right)^T M^{-1}_{D,\lambda}B_\Delta {\bf w}_\Delta\\
& = & \frac{1}{h^2}\left< g_{p_\Gamma}, g_{p_\Gamma} \right>+\left(B_\Delta {\bf w}_\Delta
\right)^T B_{\Delta, D} H_{\Delta} B_{\Delta, D}^T B_\Delta {\bf w}_\Delta\\
& = & \frac{1}{h^2}\left< g_{p_\Gamma}, g_{p_\Gamma} \right>+\left<P_{D,2} w,P_{D,2} w\right>_{\Atilde} \le \frac{1}{h^2}\left< g_{p_\Gamma}, g_{p_\Gamma} \right>+\Phi_2(H,h)\left<w,w\right>_{\Atilde},
\end{eqnarray*}
where we have used Lemma \LA{jumptwo} for the last inequality. It is  sufficient to bound the first term of the right-hand side in the above inequality.

We denote $\vvec{w} = \left( {\bf w}_I, ~{\bf w}_{\Delta}, ~{\bf w}_{\Pi} \right) \in \vvec{\Wtilde}$. Since $B_{I I} \vvec{w}_I + B_{I\Delta} \vvec{w}_\Delta + B_{I\Pi} \vvec{w}_\Pi = 0$, we have
\[
\left<g_{p_\Gamma}, g_{p_\Gamma}\right> = \left[ \begin{array}{c} B_{I I} \vvec{w}_I + B_{I\Delta} \vvec{w}_\Delta + B_{I\Pi} \vvec{w}_\Pi \\ B_{\Gamma I} {\bf w}_I + B_{\Gamma \Delta} {\bf w}_\Delta + B_{\Gamma\Pi} {\bf w}_\Pi \end{array} \right]^T
\left[ \begin{array}{c} B_{I I} \vvec{w}_I + B_{I\Delta} \vvec{w}_\Delta + B_{I\Pi} \vvec{w}_\Pi \\ B_{\Gamma I} {\bf w}_I + B_{\Gamma \Delta} {\bf w}_\Delta + B_{\Gamma\Pi} {\bf w}_\Pi \end{array} \right] = \big<\Btilde {\bf w},  \Btilde {\bf w} \big>,
\]
where  $\Btilde$ is defined in \EQ{Btilde}. From \EQ{massmatrix} and the stability of $\Btilde$, cf. Lemma \LA{BtildeStability}, we have
\begin{eqnarray}
\frac{1}{h^2}\left<g_{p_\Gamma}, g_{p_\Gamma}\right> & = & \frac{1}{h^2} \left<\Btilde {\bf w},  \Btilde {\bf w} \right> \leq C \left<\Btilde {\bf w},  \Btilde {\bf w} \right>_{Z^{-1}} =  C \max_{q \in Q} \frac{\left<\Btilde {\bf w}, q \right>^2}{\left<q, q\right>_Z} \label{equation:boundBu}\\
& \le & C \max_{q \in Q} \frac{| \vvec{w} |^2_{H^1}  \| q \|^2_{L^2}}{\| q \|^2_{L^2}} = C | \vvec{w} |^2_{H^1}
= C \left<w,w\right>_{\Atilde}, \nonumber
\end{eqnarray}
where for the last equality, we used the fact that $B_{I I} \vvec{w}_I
+ B_{I\Delta} \vvec{w}_\Delta + B_{I\Pi} \vvec{w}_\Pi = 0$ and \EQ{wg0}.
\endproof

\begin{mylemma}
\label{lemma:lowerGDiri} Let the coarse level primal velocity space ${\bf W}_\Pi$ be chosen such that \EQ{divergencefree} is enforced. For any given $y = (g_{p_{\Gamma}}, g_\lambda) \in R_G$, there exits
$w \in \Wtilde_{0}$, such that $B_C w = y$, and $\left <\Atilde w, w\right> \le \frac{C}{\beta^2}  \left< M_{D}^{-1}y, y \right>$.
\end{mylemma}

\beginproof
Given $y = (g_{p_{\Gamma}}, g_\lambda) \in R_G$, let ${\u}_{\Delta}^{(I)} = B_{\Delta, D}^T g_\lambda$ and ${\bf u}_{\Pi}^{(I)}={\bf 0}$. On each subdomain $\Omega_i$, denote $\left({\bf u}_{I}^{(I, i)}, ~p_{I}^{(I, i)}, ~p_{\Gamma}^{(I,i)}\right)$ the part obtained through the solution of~\EQ{schurplus} with given subdomain boundary values ${\u}_{\Delta}^{(i)} ={\u}_{\Delta}^{(I, i)}$ and ${\u}_{\Pi}^{(i)} = {\bf 0}$. Let $w^{(I,i)} = \left(\u^{(I,i)}_I, ~p_I^{(I,i)}, ~{\bf u}_{\Delta}^{(I,i)}, ~ {\bf u}^{(I,i)}_\Pi\right)$, the corresponding global vectors
$w^{(I)} = \left(\u^{(I)}_I, ~p_I^{(I)}, ~{\bf u}_{\Delta}^{(I)}, ~ {\bf u}^{(I)}_\Pi\right)$, and $\u^{(I)} = \left(\u^{(I)}_I,  ~{\bf u}_{\Delta}^{(I)}, ~ {\bf u}^{(I)}_\Pi\right)$. Then we know from \EQ{schurplus} that $w^{(I)}\in \widetilde{W}_{0}$,
and
\be
\label{equation:bcID1}
B_C w^{(I)} = \left[
\begin{array}{cccc}
B_{\Gamma I} & 0 & B_{\Gamma \Delta} & B_{\Gamma \Pi} \\[0.8ex]
0            & 0 & B_{\Delta}        & 0              \end{array}
\right] \left[ \begin{array}{c}
{\bf u}_I^{(I)}        \\[0.8ex]
p_I^{(I)}              \\[0.8ex]
{\bf u}_{\Delta}^{(I)} \\[0.8ex]
{\bf u}_{\Pi}^{(I)}
\end{array} \right] = \left[ \begin{array}{c}
0     \\[0.8ex] g_\lambda \end{array} \right],
\ee
where we used the fact that $B_\Delta B_{\Delta, D}^T = I$.

Using Lemma \LA{stokesnormPlus}, we have
\[
| \u^{(I, i)} |^2_{H^1(\Omega_i)} = \left| \left[ \begin{array}{c}
     {\bf u}_{\Delta}^{(I, i)} \\ 0 \end{array} \right]
\right|^2_{S_{\Delta\Pi}^{(i)}} \leq \frac{C}{\beta^2} \left|
 \left[ \begin{array}{c} {\bf u}_{\Delta}^{(I, i)} \\ 0 \end{array}
 \right] \right|^2_{H^{1/2}(\partial \Omega^i)},
\]
and summing over the subdomains,
\be \label{equation:Iprop}
| \u^{(I)} |^2_{H^1} \leq
\frac{C}{\beta^2} \sum_{i=1}^N \left|
 \left[ \begin{array}{c} {\bf u}_{\Delta}^{(I, i)} \\ 0 \end{array}
 \right] \right|^2_{H^{1/2}(\partial \Omega^i)}.
\ee

We consider a solution to the following fully assembled system of linear equations of the form~\EQ{matrix}: find
$\left({\bf u}_I^{(II)}, ~p_I^{(II)}, ~{\bf u}_{\Gamma}^{(II)}, ~p_\Gamma^{(II)}\right) \in {\bf W}_I \bigoplus Q_I \bigoplus {\bf W}_{\Gamma} \bigoplus Q_\Gamma$, such that
\be
\label{equation:fassembled}
\left[
\begin{array}{cccc}
A_{II}      & B_{II}^T       & A_{I \Gamma}      & B_{\Gamma I}^T     \\[0.8ex]
B_{II}      & 0              & B_{I \Gamma}      & 0                  \\[0.8ex]
A_{\Gamma I}& B_{I \Gamma} ^T& A_{\Gamma\Gamma}  & B_{\Gamma \Gamma}^T \\[0.8ex]
B_{\Gamma I}& 0              & B_{\Gamma \Gamma} & 0
\end{array}
\right]
\left[ \begin{array}{c}
{\bf u}_I^{(II)}        \\[0.8ex]
p_I^{(II)}              \\[0.8ex]
{\bf u}_{\Gamma}^{(II)} \\[0.8ex]
p_{\Gamma}^{(II)}
\end{array} \right] =
\left[ \begin{array}{l}
{\bf 0}        \\[0.8ex]
0             \\[0.8ex]
{\bf 0}        \\[0.8ex]
g_{p_{\Gamma}}
\end{array} \right] \mbox{ , }
\ee
where a particular right-hand side is chosen.  Since $y\in R_G$, we can see from \EQ{Grange} and \EQ{dfree}, that $g_{p_\Gamma}^T 1_{{p_\Gamma}} =0$. Therefore the above system has a solution.

Denote ${\bf u}^{(II)} = \left( {\bf u}_I^{(II)}, ~{\bf u}_{\Gamma}^{(II)} \right)$. Then from Lemma \ref{lemma:paul} and \EQ{massmatrix}, we have
\be \label{equation:uIboundD}
| {\bf u}^{(II)} |^2_{H^1} \leq \frac{1}{\beta^2} \left\| \left[ \begin{array}{l} 0 \\[0.8ex]
g_{p_\Gamma} \end{array} \right] \right\|^2_{Z^{-1}} \le \frac{C}{\beta^2 h^2}\left<g_{p_\Gamma},g_{p_\Gamma}\right>.
\ee

Split the continuous subdomain boundary velocity ${\bf u}_{\Gamma}^{(II)}$ into the dual part ${\bf u}_{\Delta}^{(II)}$ and the primal part ${\bf u}_{\Pi}^{(II)}$, and denote $w^{(II)} = \left({\bf u}_I^{(II)}, ~p_I^{(II)}, ~{\bf u}_{\Delta}^{(II)}, ~{\bf u}_{\Pi}^{(II)}\right)$. Then we have from \EQ{fassembled} that $w^{(II)}\in \widetilde{W}_{0}$ and
\be
\label{equation:bcID2}
B_C w^{(II)} = \left[
\begin{array}{cccc}
B_{\Gamma I} & 0 & B_{\Gamma \Delta} & B_{\Gamma \Pi} \\[0.8ex]
0            & 0 & B_{\Delta}        & 0              \end{array}
\right] \left[ \begin{array}{c}
{\bf u}_I^{(II)}        \\[0.8ex]
p_I^{(II)}              \\[0.8ex]
{\bf u}_{\Delta}^{(II)} \\[0.8ex]
{\bf u}_{\Pi}^{(II)}
\end{array} \right] = \left[ \begin{array}{c}
g_{p_\Gamma}      \\[0.8ex] 0  \end{array} \right].
\ee

Let $w = w^{(I)} + w^{(II)}$. We have $w \in \Wtilde_{0}$, $B_C w = y$, from \EQ{bcID1} and \EQ{bcID2},  and from \EQ{wg0}
\[
| w |^2_{\widetilde{A}}=| \u^{(I)}+ \u^{(II)}|^2_{H^1}  \le| \u^{(I)} |^2_{H^1}+| \u^{(II)} |^2_{H^1} \leq
\frac{C}{\beta^2} \sum_{i=1}^N \left|
 \left[ \begin{array}{c} {\bf u}_{\Delta}^{(I, i)} \\ 0 \end{array}
 \right] \right|^2_{H^{1/2}(\partial \Omega^i)} + \frac{C}{\beta^2 h^2}\left<g_{p_\Gamma},g_{p_\Gamma}\right>,
\]
where we have used \EQ{Iprop} and \EQ{uIboundD} in the last inequality.

On the other hand, we have from \EQ{subdomainwise}
\begin{eqnarray*}
\left< M_{D}^{-1}y,y \right>&=&
\frac{1}{h^2}\left<g_{p_\Gamma}, g_{p_\Gamma}\right>+ g_\lambda^T
M^{-1}_{D,\lambda} g_\lambda\\
& = & \frac{1}{h^2}\left<g_{p_\Gamma}, g_{p_\Gamma}\right>+ g_\lambda^T \left( \sum_{i=1}^N  
B_{\Delta, D}^{(i)}H_{\Delta}^{(i)} B_{\Delta, D}^{(i)^T}
\right) g_\lambda \\
&=&\frac{1}{h^2}\left<g_{p_\Gamma},g_{p_\Gamma}\right>+\sum_{i=1}^N \left|
 \left[ \begin{array}{c} {\bf u}_{\Delta}^{(I,i)} \\ 0 \end{array}
 \right]
\right|^2_{H^{1/2}(\partial \Omega^i)}. \qquad \Box
\end{eqnarray*}

With Lemmas \ref{lemma:upperGDiri} and \ref{lemma:lowerGDiri}, similar to the proof of Theorem~\ref{theorem:tcondF}, we have the following theorem.
\begin{mytheorem}
\label{theorem:DiricondG} Let the coarse level primal velocity space ${\bf W}_\Pi$ be chosen such that \EQ{divergencefree} is enforced. For all $x = (p_{\Gamma}, ~\lambda) \in R_{M_{D}^{-1} G}$,
\[
c \beta^2 \left<M_{D}x,x \right>\leq \left< G x,x \right> \leq C \Phi_2(H,h) \left<
M_{D}x,x \right>,
\]
where $\Phi_2(H,h)$ is as defined in Lemma \LA{jumptwo}, $\beta$ is the inf-sup constant specified in \EQ{infsup}.
\end{mytheorem}

When
$p_\Gamma$ in \EQ{spdG} is empty, the matrix G contains only the second diagonal block in \EQ{Gtwo} and the Dirichlet preconditioner for \EQ{spdG} becomes only $M^{-1}_{D,\lambda}$, as given in \EQ{subdomainwise}, i.e.,
\[
M^{-1}_{D} = M^{-1}_{D,\lambda} = B_{\Delta, D} H_{\Delta}  B_{\Delta, D}^T.
\]
Counterparts of Lemmas \LA{upperGDiri} and \LA{lowerGDiri} for this case can be proved as well; in fact, their proofs are essentially the first half in the proofs of Lemmas \LA{upperGDiri} and \LA{lowerGDiri}, respectively. The condition number bound in Theorem \ref{theorem:DiricondG} can then be established for this case in the same way.

\section{Numerical experiments}
\label{section:numerics}

We consider solving the incompressible Stokes problem \EQ{Stokes} in the  domain $\Omega=[0,1]\times
[0,1]$. Zero Dirichlet boundary condition is used. The right-hand side function $\vvec{f}$ is chosen such that the exact solution is
$$\u=\left[\begin{array}{c}
\sin^3(\pi x)\sin^2(\pi y)\cos(\pi y)\\[0.8ex]
-\sin^2(\pi x)\sin^3(\pi y)\cos(\pi x)
\end{array}\right] \quad \mbox{and}\quad
p=x^2-y^2.
$$

Two mixed finite element discretizations, as shown on Figures \ref{fig:2dfem} and \ref{figure:TaylorHood} in Section \ref{section:FEM}, are used for the cases of discontinuous and continuous pressures, respectively. The preconditioned system \EQ{Mspd} is solved by the preconditioned conjugate gradient method; the iteration is stopped when the $L^2-$norm of the residual is reduced by a factor of $10^{-6}$.

In each of the following tables, we present the performance of three
different variants of the FETI-DP algorithm represented under the same
framework, as discussed in Sections~\ref{subsection:cp} and~\ref{subsection:dp}: ``continuous pressure" for the case when the
continuous pressure is used in the algorithm and $p_\Gamma$ contains
all the subdomain boundary pressure degrees of freedom;
``discontinuous pressure" for the case when the discontinuous pressure
is used and $p_\Gamma$ contains just one pressure degree of freedom
from each subdomain; ``$p_\Gamma$ empty" for the case when the
discontinuous pressure is used and $p_\Gamma$ is chosen empty. For
each case, the extreme eigenvalues and the iteration count for each
experiment are shown. The two methods discussed in \cite{li05}
and~\cite{kim10} solve the same system \EQ{spdG},  and their only difference is in
the  implementation of multiplying $G$ by a vector, cf. Remark
\ref{remark:method12}. Therefore  their convergence rates are the same,
when they are equipped with the same type preconditioner, and their performance is reported under the case ``$p_\Gamma$ empty" in the tables.

Tables \ref{table:One} and \ref{table:Two} show the performance of
using the lumped preconditioner for different cases with two choices
of the coarse level prime variables.  In Table \ref{table:One}, the
first choice is used, namely only the subdomain corner velocities are
taken as the coarse level primal variables. We can see that for each
variant of the FETI-DP algorithm, the convergence rate is independent
of the number of subdomains for fixed $H/h$; for fixed number of
subdomains, the condition number grows presumably in the order of
$(H/h)\left( 1 + \log(H/h)\right)$ as established in Section
\ref{section:lumped}. In Table~\ref{table:Two}, we test the lumped
preconditioner with the second choice of  the coarse level primal
space, namely it  contains both the subdomain corner velocity variables and the edge-average velocity variables such that \EQ{divergencefree} is enforced, as discussed in Section \ref{section:DDM}. Even though the edge-average velocity variables are not required in the coarse level primal space for the lumped preconditioner case, including them improves the convergence rate for each method. We also observe from Tables \ref{table:One} and \ref{table:Two} that performances of the three variants of the FETI-DP algorithm are quite similar, while the convergence when using discontinuous pressure and choosing $p_\Gamma$ empty is a little faster than the other two cases.

Tables \ref{table:Three} and \ref{table:Four} show the performance of
using the Dirichlet preconditioner. In Table \ref{table:Three}, only
the subdomain corner velocities are taken as the coarse level primal
variables, for which the divergence free boundary condition
\EQ{divergencefree} is not satisfied and no scalable condition number
bound of the FETI-DP algorithm is available. Indeed Table
\ref{table:Three} shows that for each variant of the FETI-DP
algorithm, the convergence rate deteriorates with the increase of the
number of subdomains when the subdomain problem size is fixed. In Table \ref{table:Four}, both the subdomain corner and the edge-average velocity variables are taken as the coarse level primal variables such that \EQ{divergencefree} is enforced. We can see that for each variant of the FETI-DP algorithm, the convergence rate is independent of the number of subdomains for fixed $H/h$; for fixed number of subdomains, the condition number grows presumably in the order of $\left( 1 + \log(H/h)\right)^2$ as established in Section~\ref{section:Dirichlet}.

\begin{table}[t]
\caption{\label{table:One} Performance using lumped preconditioner $M^{-1}_{L}$, with corner primal variables.}
\centering
\begin{tabular}{ccccccccccccc}
\hline
&&\multicolumn{3}{c}{continuous pressure} & & \multicolumn{3}{c}{discontinuous pressure} & &\multicolumn{3}{c}{$p_\Gamma$ empty}\\
\cline{3-5} \cline{7-9} \cline{11-13}
$H/h$ & \#sub & $\lambda_{min}$ & $\lambda_{max}$ & iter & & $\lambda_{min}$
& $\lambda_{max}$ & iter  & & $\lambda_{min}$
& $\lambda_{max}$ & iter \\[0.8ex]
\hline \\
8& $4  \times  4$ & 0.35 &  8.92 & 21 & & 0.48 & 7.93 & 22 &  & 0.56 & 7.37 & 20 \\[0.8ex]
 & $8  \times  8$ & 0.35 & 10.07 & 28 & & 0.48 & 9.00 & 25 &  & 0.56 & 8.46 & 22 \\[0.8ex]
 & $16 \times 16$ & 0.35 & 10.23 & 29 & & 0.48 & 9.20 & 25 &  & 0.56 & 8.71 & 22 \\[0.8ex]
 & $24 \times 24$ & 0.35 & 10.30 & 29 & & 0.48 & 9.20 & 25 &  & 0.56 & 8.68 & 22 \\[0.8ex]
 & $32 \times 32$ & 0.35 & 10.33 & 29 & & 0.48 & 9.21 & 25 &  & 0.56 & 8.68 & 22 \\
\hline \\
\#sub & $H/h$ &  $\lambda_{min}$ & $\lambda_{max}$ & iter & &  $\lambda_{min}$ & $\lambda_{max}$ & iter & & $\lambda_{min}$
& $\lambda_{max}$ & iter \\[0.8ex]
\hline \\
$8 \times  8$&  4 & 0.30 &  4.22 & 21 & & 0.41 &  3.91 & 19 &  & 0.54 &  3.75 & 16 \\[0.8ex]
             &  8 & 0.35 & 10.07 & 28 & & 0.48 &  9.00 & 25 &  & 0.56 &  8.46 & 22 \\[0.8ex]
             & 16 & 0.35 & 24.22 & 36 & & 0.49 & 21.39 & 36 &  & 0.56 & 20.30 & 33 \\[0.8ex]
             & 24 & 0.35 & 40.12 & 43 & & 0.50 & 35.56 & 43 &  & 0.57 & 33.89 & 39 \\[0.8ex]
             & 32 & 0.35 & 57.15 & 50 & & 0.50 & 50.87 & 50 &  & 0.57 & 48.62 & 45 \\
\hline
\end{tabular}
\end{table}

\begin{table}[t]
\caption{\label{table:Two} Performance using lumped preconditioner $M^{-1}_{L}$, with corner and edge-average primal variables.}
\centering
\begin{tabular}{ccccccccccccc}
\hline
&&\multicolumn{3}{c}{continuous pressure} & & \multicolumn{3}{c}{discontinuous pressure} & &\multicolumn{3}{c}{$p_\Gamma$ empty}\\
\cline{3-5} \cline{7-9} \cline{11-13}
$H/h$ & \#sub & $\lambda_{min}$ & $\lambda_{max}$ & iter & & $\lambda_{min}$
& $\lambda_{max}$ & iter  & & $\lambda_{min}$
& $\lambda_{max}$ & iter \\[0.8ex]
\hline \\
8& $4  \times  4$ & 0.36 & 4.29 & 17  & & 0.48 & 3.78 & 17 & & 0.56 & 3.39 & 14 \\[0.8ex]
 & $8  \times  8$ & 0.36 & 5.29 & 21  & & 0.49 & 4.47 & 18 & & 0.56 & 4.01 & 16 \\[0.8ex]
 & $16 \times 16$ & 0.36 & 5.56 & 21  & & 0.49 & 4.68 & 19 & & 0.56 & 4.29 & 16 \\[0.8ex]
 & $24 \times 24$ & 0.36 & 5.61 & 21  & & 0.50 & 4.77 & 19 & & 0.55 & 4.42 & 16 \\[0.8ex]
 & $32 \times 32$ & 0.36 & 5.64 & 21  & & 0.50 & 4.80 & 19 & & 0.55 & 4.46 & 16 \\
\hline \\
\#sub & $H/h$ &  $\lambda_{min}$ & $\lambda_{max}$ & iter & &  $\lambda_{min}$ & $\lambda_{max}$ & iter & & $\lambda_{min}$
& $\lambda_{max}$ & iter \\[0.8ex]
\hline \\
$8  \times 8$&  4 & 0.33 &  4.00 & 18 & & 0.43 &  2.80 & 16 & & 0.55 &  1.91 & 11 \\[0.8ex]
             &  8 & 0.36 &  5.29 & 21 & & 0.49 &  4.47 & 18 & & 0.56 &  4.01 & 16 \\[0.8ex]
             & 16 & 0.36 & 11.63 & 26 & & 0.50 &  9.85 & 26 & & 0.56 &  9.31 & 23 \\[0.8ex]
             & 24 & 0.36 & 18.67 & 31 & & 0.50 & 16.05 & 32 & & 0.57 & 15.36 & 29 \\[0.8ex]
             & 32 & 0.36 & 26.12 & 36 & & 0.50 & 22.67 & 37 & & 0.57 & 21.83 & 33 \\
\hline
\end{tabular}
\end{table}

\begin{table}[t]
\caption{\label{table:Three}  Performance using Dirichlet preconditioner $M^{-1}_{D}$, with corner primal variables.}
\centering
\begin{tabular}{ccccccccccccc}
\hline
&&\multicolumn{3}{c}{continuous pressure} & & \multicolumn{3}{c}{discontinuous pressure} & &\multicolumn{3}{c}{$p_\Gamma$ empty}\\
\cline{3-5} \cline{7-9} \cline{11-13}
$H/h$ & \#sub & $\lambda_{min}$ & $\lambda_{max}$ & iter & & $\lambda_{min}$
& $\lambda_{max}$ & iter  & & $\lambda_{min}$
& $\lambda_{max}$ & iter \\[0.8ex]
\hline \\
8& $4  \times  4$ & 0.32 & 3.11 & 18 & & 0.38 & 2.83 & 16 & & 0.42 & 1.87  & 12 \\[0.8ex]
 & $8  \times  8$ & 0.29 & 3.42 & 19 & & 0.33 & 3.01 & 18 & & 0.33 & 2.19  & 14 \\[0.8ex]
 & $16 \times 16$ & 0.25 & 3.52 & 21 & & 0.28 & 3.09 & 19 & & 0.28 & 2.30  & 16 \\[0.8ex]
 & $24 \times 24$ & 0.24 & 3.56 & 22 & & 0.26 & 3.11 & 19 & & 0.26 & 2.35  & 16 \\[0.8ex]
 & $32 \times 32$ & 0.23 & 3.57 & 22 & & 0.25 & 3.12 & 20 & & 0.25 & 2.37  & 17 \\
\hline \\
\#sub & $H/h$ &  $\lambda_{min}$ & $\lambda_{max}$ & iter & &  $\lambda_{min}$ & $\lambda_{max}$ & iter & & $\lambda_{min}$
& $\lambda_{max}$ & iter \\[0.8ex]
\hline \\
$8 \times  8$&  4 & 0.27 & 3.82 & 22  & & 0.33 & 2.73 & 18 &  & 0.37 & 1.70 & 13 \\[0.8ex]
             &  8 & 0.29 & 3.42 & 19  & & 0.33 & 3.01 & 18 &  & 0.33 & 2.19 & 14 \\[0.8ex]
             & 16 & 0.30 & 4.00 & 21  & & 0.32 & 3.39 & 19 &  & 0.32 & 2.84 & 16 \\[0.8ex]
             & 24 & 0.30 & 4.39 & 22  & & 0.32 & 3.69 & 19 &  & 0.32 & 3.28 & 17 \\[0.8ex]
             & 32 & 0.31 & 4.71 & 23  & & 0.32 & 3.95 & 20 &  & 0.32 & 3.62 & 19 \\
\hline
\end{tabular}
\end{table}

\begin{table}[t]
\caption{\label{table:Four} Performance using Dirichlet preconditioner $M^{-1}_{D}$, with corner and edge-average primal variables.}
\centering
\begin{tabular}{ccccccccccccc}
\hline
&&\multicolumn{3}{c}{continuous pressure} & & \multicolumn{3}{c}{discontinuous pressure} & &\multicolumn{3}{c}{$p_\Gamma$ empty}\\
\cline{3-5} \cline{7-9} \cline{11-13}
$H/h$ & \#sub & $\lambda_{min}$ & $\lambda_{max}$ & iter & & $\lambda_{min}$
& $\lambda_{max}$ & iter  & & $\lambda_{min}$
& $\lambda_{max}$ & iter \\[0.8ex]
\hline \\
8& $4  \times  4$ & 0.35 & 2.84 & 16 & & 0.40 & 2.53 & 15 & & 0.47 & 1.35 & 10 \\[0.8ex]
 & $8  \times  8$ & 0.35 & 2.96 & 16 & & 0.41 & 2.65 & 15 & & 0.47 & 1.59 & 10 \\[0.8ex]
 & $16 \times 16$ & 0.35 & 3.00 & 15 & & 0.42 & 2.71 & 15 & & 0.47 & 1.75 & 11 \\[0.8ex]
 & $24 \times 24$ & 0.35 & 3.02 & 15 & & 0.42 & 2.74 & 15 & & 0.47 & 1.79 & 11 \\[0.8ex]
 & $32 \times 32$ & 0.35 & 3.03 & 15 & & 0.43 & 2.75 & 15 & & 0.47 & 1.82 & 11 \\
\hline \\
\#sub & $H/h$ &  $\lambda_{min}$ & $\lambda_{max}$ & iter & &  $\lambda_{min}$ & $\lambda_{max}$ & iter & & $\lambda_{min}$
& $\lambda_{max}$ & iter \\[0.8ex]
\hline \\
$8 \times  8$&  4 & 0.33 & 3.67 & 18 & & 0.41 & 2.60 & 16 & & 0.46 & 1.30  & 10 \\[0.8ex]
             &  8 & 0.35 & 2.96 & 16 & & 0.41 & 2.65 & 15 & & 0.47 & 1.59  & 10 \\[0.8ex]
             & 16 & 0.35 & 2.87 & 15 & & 0.40 & 2.81 & 15 & & 0.47 & 2.00  & 11 \\[0.8ex]
             & 24 & 0.34 & 2.99 & 16 & & 0.40 & 2.97 & 16 & & 0.50 & 2.30  & 12 \\[0.8ex]
             & 32 & 0.34 & 3.21 & 16 & & 0.40 & 3.11 & 16 & & 0.48 & 2.52  & 13 \\
\hline
\end{tabular}
\end{table}

\end{document}